\documentclass[12pt]{amsart}
\usepackage{amsthm}
\usepackage{amsmath, mathrsfs}
\usepackage{amssymb}
\usepackage{hyperref}
\usepackage{enumerate}
\usepackage{latexsym,array}
\usepackage{amsfonts}
\usepackage{shadow}
\usepackage{tikz}
\usepackage{rotating}
\usepackage{tabularx}
\usetikzlibrary{arrows}
\usepackage[english]{babel}
\usepackage{graphicx}
\usepackage[section]{placeins}
\usepackage{fullpage}

\newcommand{\bigH}{\mathcal{H}}

\newcommand{\bigF}{\mathcal{F}}

\newcommand{\bigM}{\mathcal{M}}

\newcommand{\bigO}{\mathcal{O}}
\newcommand{\ball}{\mathbb{B}_d}
\newcommand{\Cball}{\overline{\mathbb{B}_d}}

\newtheorem{Pa}{Paper}[section]
\newtheorem{Tm}[Pa]{{\bf Theorem}}
\newtheorem{TmA}{\bf Theorem}

\newtheorem{La}[Pa]{{\bf Lemma}}
\newtheorem{Cy}[Pa]{{\bf Corollary}}

\newtheorem{Pn}[Pa]{{\bf Proposition}}

\theoremstyle{definition}

\newtheorem{Ex}[Pa]{{\bf Example}}
\theoremstyle{remark}
\newtheorem{Rk}[Pa]{{\bf Remark}}
\tikzset{node distance=2cm, auto}

\newcommand{\cA}{\mathcal{A}}
\newcommand{\cD}{\mathcal{D}}
\newcommand{\cE}{\mathcal{E}}
\newcommand{\cF}{\mathcal{F}}
\newcommand{\cH}{\mathcal{H}}
\newcommand{\cI}{\mathcal{I}}
\newcommand{\cM}{\mathcal{M}}

\newcommand{\cZ}{\mathcal{Z}}

\newcommand{\bB}{\mathbb{B}}
\newcommand{\bC}{\mathbb{C}}

\newcommand{\bN}{\mathbb{N}}

\newcommand{\ol}{\overline}
\newcommand{\spn}{\textrm{span}}
\newcommand{\alg}{\textrm{Alg}}
\newcommand{\mlt}{\textrm{Mult}}
\newcommand{\Aut}{\textrm{Aut}}

\newcommand{\be}{\begin{equation}}
\newcommand{\ee}{\end{equation}}

\date{}

\title[The isomorphism problem for complete Pick algebras]
{The isomorphism problem for complete Pick algebras:\\ a survey}

\author[G. Salomon]{Guy Salomon}
\address{Department of Mathematics\\Technion - Israel Institute of Technology\\Haifa 3200003\\Israel}
\email{guy.salomon@tx.technion.ac.il}

\author[O.M. Shalit]{Orr Moshe Shalit}

\address{Department of Mathematics\\Technion - Israel Institute of Technology\\Haifa 3200003\\Israel}
\email{oshalit@tx.technion.ac.il}

\thanks{The second author was partially supported by ISF Grant no. 474/12, by
EU FP7/2007-2013 Grant no. 321749, and by GIF Grant no. 2297-2282.6/20.1.}

\begin{document}
\maketitle


\begin{abstract}
Complete Pick algebras --- these are, roughly, the multiplier algebras in which Pick's interpolation theorem holds true --- have been the focus of much research in the last twenty years or so.
All (irreducible) complete Pick algebras may be realized concretely as the algebras obtained by restricting multipliers on Drury-Arveson space to a subvariety of the unit ball;
to be precise: every irreducible complete Pick algebra has the form $\cM_V = \{f\big|_V : f \in \cM_d\}$, where $\cM_d$ denotes the multiplier algebra of the Drury-Arveson space $H^2_d$, and $V$ is the joint zero set of some functions in $\cM_d$.
In recent years several works were devoted to the classification of complete Pick algebras in terms of the complex geometry of the varieties with which they are associated.
The purpose of this survey is to give an account of this research in a comprehensive and unified way.
We describe the array of tools and methods that were developed for this program, and take the opportunity to clarify, improve, and correct some parts of the literature.
\end{abstract}

\section{Introduction}

\subsection{Motivation and background}

Consider the following two classical theorems.
\begin{TmA}[Gelfand, \cite{GRS}]\label{thm:A}
Let $X$ and $Y$ be two compact Hausdorff spaces. The algebras of continuous functions $C(X)$ and $C(Y)$ are isomorphic if and only if $X$ and $Y$ are homeomorphic.
\end{TmA}

\begin{TmA}[Bers, \cite{Bers}]\label{thm:B}
Let $U$ and $V$ be open subsets of $\bC$. The algebras of holomorphic functions ${\rm Hol}(U)$ and ${\rm Hol}(V)$ are isomorphic if and only if $U$ and $V$ are biholomorphic
\end{TmA}

The common theme of these two theorems is that an appropriate algebra of functions on a space encapsulates in its algebraic structure every aspect of the topological/complex-geometric structure of the space. The problem that we are concerned with in this paper has a very similar flavour.
Let $\cM_d$ denote the algebra of multipliers on Drury-Arveson space --- precise definitions will be given in the next section, for now it suffices to say that $\cM_d$ is a certain algebra of bounded analytic functions on the unit ball $\bB_d \subseteq \bC^d$.
For every analytic variety $V \subseteq \bB_d$ one may define the algebra
\[
\cM_V = \{ f\big|_V : f \in \cM_d\}.
\]
The natural question to ask is: in what ways does the variety $V$ determine the algebra $\cM_V$, and vice versa?
In other words, if $\cM_V$ and $\cM_W$ are algebraically isomorphic, can we conclude that $V$ and $W$ are ``isomorphic" in some sense?
Conversely, if $V$ and $W$ are, say, biholomorphic, can we conclude that the algebras are isomorphic?

As we shall explain below, $\cM_V$ is also an operator algebra: it is the multiplier algebra of a certain reproducing kernel Hilbert space on $V$, and it is generated by the multiplication operators $[M_{z_i}  h](z) =   z_i h(z)$ (it will be convenient to denote henceforth $Z_i = M_{z_i}$).
Thus one can ask: do the Banach algebraic or operator algebraic structures of $\cM_V$  encode finer complex-geometric aspects of $V$?

These questions in themselves are interesting, natural, nontrivial, and studying them involves a collection of tools combining function theory, complex geometry and operator theory.
However, it is worth noting that there are routes, other than analogy with Theorems \ref{thm:A} and \ref{thm:B}, that lead one to study the structure and classify the algebras $\cM_V$ described above.

One path that leads to considering the algebras $\cM_V$ comes from non-selfadjoint operator algebras: it is the study of operator algebras universal with respect to some polynomial relations.
For simplicity consider the case in which $V = \cZ_{\bB_d}(\cI)$ is the zero set of a radical and homogeneous polynomial ideal $\cI \triangleleft \bC[z_1, \ldots, z_d]$, where
\[
\cZ_{\bB_d}(\cI) = \{\lambda \in \bB_d \mid p(\lambda) = 0 \,\textrm{ for all }\, p \in \cI\}.
\]
Then $\cM_V$ is the universal \textsc{wot}-closed unital operator algebra, that is generated by a pure commuting row contraction $T = (T_1, \ldots, T_d)$ satisfying the relations in $\cI$ (see \cite{Pop06, ShalitSolel}).
This means that
\begin{enumerate}
\item The $d$-tuple of operators $\left(Z_1, \ldots, Z_d\right)$, given by multiplication by the coordinate functions, is a pure, commuting row contraction satisfying the relations in $\cI$, and it generates $\cM_V$;
\item For any such tuple $T$, there is a unital, completely contractive and \textsc{wot}-continuous homomorphism from $\cM_V$ into $\overline{\alg}^{\textsc{wot}}\left(1,T\right)$ determined by $Z_i \mapsto T_i$.
\end{enumerate}
In general (when $V$ is not necessarily the variety of a homogeneous polynomial ideal) it is a little more complicated to explain the universal property of $\cM_V$. Roughly, $\cM_V$ is universal for tuples ``satisfying the relations'' in $J_V = \{f \in \cM_d \mid f(\lambda) = 0 \,\textrm{ for all } \, \lambda \in V\}$.

Thus the algebras $\cM_V$ are an operator algebraic version of the coordinate ring on an algebraic variety, and studying the relations between the structure of $\cM_V$ and the geometry of $V$ can be considered as rudimentary steps in developing ``operator algebraic geometry''.

A different road that leads one to consider the collection of algebras $\cM_V$ runs from function theory, in particular from the theory of Pick interpolation.
Let $\cH$ be a reproducing kernel Hilbert space on a set $X$ with kernel $k$.
If $x_1, \ldots, x_n \in X$ and $A_1, \ldots, A_n \in M_k(\bC)$, then one may consider the problem of finding a matrix valued multiplier $F : X \rightarrow M_k(\bC)$ which has multiplier norm $1$ and satisfies
\[
F(x_i) = A_i \, \, , \,\, i = 1 , \ldots, d.
\]
This is called the {\em Pick interpolation problem}.
It is not hard to show that a necessary condition for the existence of such a multiplier is that the following matrix inequality hold:
\be\label{eq:pick}
\left[(1 - F(x_i) F(x_j)^*) K(x_i,x_j)\right]_{i,j=1}^n \geq 0 .
\ee
G. Pick showed that for the Szeg\H{o} kernel $k(z,w) = (1-z \bar w)^{-1}$ the condition \eqref{eq:pick} is also a sufficient condition for the existence of a solution to this problem \cite{Pick}.
Kernels for which condition \eqref{eq:pick} is a sufficient condition for the existence of a solution to the Pick interpolation problem have come to be called {\em complete Pick kernels}, and their multiplier algebras {\em complete Pick algebras}.
We refer the reader to the monograph \cite{AglerMcCarthy_PickInterpolation} for thorough introduction to Pick interpolation and complete Pick kernels.
The connection to our problem is the following theorem, which states that under a harmless irreducibility assumption all complete Pick algebras are completely isometrically isomorphic to one of the algebras $\cM_V$ described above.

\begin{TmA}[Agler-M\raise.45ex\hbox{c}Carthy, \cite{AM00}]\label{thm:AM}
Let $\cH$ be a reproducing kernel Hilbert space with an irreducible complete kernel $k$.
Then there exists $d \in \bN \cup \{\infty\}$ and there is an analytic subvariety $V  \subseteq \bB_d$ such that the multiplier algebra $\mlt({\cH})$ of $\cH$ is unitarily equivalent to $\cM_V$.
\end{TmA}
In fact the theorem of Agler-M\raise.45ex\hbox{c}Carthy says much more: the Hilbert space $\cH$ can (up to some rescaling) be considered as a Hilbert space of functions on $V$, which is a subspace of the Drury-Arveson space.
Since we require this result only for motivation, we do not go into further detail.

Thus, by studying the algebras $\cM_V$ in terms of the complex-geometric structure of $V$ one may hope to obtain a structure theory of irreducible complete Pick algebras.
In particular, we may hope to use the varieties as complete invariants of irreducible complete Pick algebras up to isomorphism --- be it algebraic, isometric or spatial.
This is why we call this study {\em The Isomorphism Problem for Complete Pick Algebras}.

\subsection{About this survey}

The goal of this survey is to present in a unified way the main results on the isomorphism problem for complete Pick algebras obtained in recent years.
We do not provide all the proofs, but we do give proofs (or at least an outline) to most key results, in order to highlight the techniques involved.
We give precise references so that all omitted details can be readily found by the interested reader.
We also had to omit some results, but all results directly related to this survey may be found in the cited references.

Although one may treat the case where $V \subseteq \bB_d$ and $W \subseteq \bB_{d'}$ where $d$ and $d'$ might be different, we will only treat the case where $d = d'$.
It is easy to see that this simplification results in no real loss.

This paper also contains some modest improvements to the results appearing in the literature.
In some cases we unify, in others we simplify the proof somewhat, in one case we were able to extend a result from $d<\infty$ to $d = \infty$ (see Theorem \ref{Tm:Isom=>CIS}).
There is also one case where we correct a mistake that appeared in an earlier paper (see Remark \ref{Rk:AffCodim}).

Furthermore, we take this opportunity to call to attention a little mess that resides in the literature, and try to set it right.
(The reader may skip the following paragraph and return to it after reading Section \ref{subsec:char}.)
The results we review in this survey are based directly on results in the papers \cite{APV03,ARS08,DHS_EmbeddedDiscs,DRS1,DRS2,Hartz,KerrMcCarthyShalit}.
The papers \cite{DHS_EmbeddedDiscs,DRS2} relied in a significant way on many earlier results of Davidson and Pitts \cite{DavidsonPitts_ncToeplitz,DavPittsPick,DavPitts1}, and in particular on \cite[Theorem 3.2]{DavidsonPitts_ncToeplitz}.
The content of that theorem, phrased in the language of this survey, is that over every point of $V$ there lies a unique character in the maximal ideal space $M(\cM_V)$, and moreover that there are no characters over points of $\bB_d \setminus V$.
Unfortunately, at the time that the papers \cite{DHS_EmbeddedDiscs,DRS2} were in press it was observed by Michael Hartz that  \cite[Theorem 3.2]{DavidsonPitts_ncToeplitz} is true only under the assumption $d < \infty$, a counter example shows that it is false for $d = \infty$ (see the example on the first page of \cite{DHSErr}, or Example 2.4 in the arXiv version of \cite{DHS_EmbeddedDiscs}).

Luckily, the main results of \cite{DHS_EmbeddedDiscs,DRS2} survived this disaster, but significant changes in the arguments were required, and some of the results survived in a weaker form.
The paper \cite{DHS_EmbeddedDiscs} has an erratum \cite{DHSErr}, and \cite{DRS2} contains some corrections made in proof.
However, thorough revisions of the papers \cite{DHS_EmbeddedDiscs,DRS2} appeared on the arXiv, and when we refer to these papers we refer to the arXiv versions.
We direct the interested reader to the arXiv versions.

\subsection{Overview of main results}

Sections \ref{sec:prelim} and \ref{sec:weak} contain some basic results which are used in all of the classification schemes.
The main results are presented in Sections \ref{sec:iso}, \ref{sec:alg} and \ref{sec:discs}, which can be read independently after Sections \ref{sec:prelim} and \ref{sec:weak}.
Some open problems are discussed in the final section.

The following table summarizes what is known and what is not known regarding the isomorphism problem of the algebras $\cM_V$, where $V$ is a variety in a finite dimensional ball.
(In several cases the result also holds for $d = \infty$, see caption).
\begin{sidewaystable}
\vspace{12cm}
\begin{tabular}{l  l c c l l}
\hline
\begin{minipage}[t]{0.21\linewidth}
\textbf{Conditions on $V$, $W$}
\end{minipage}
&
\begin{minipage}[t]{0.21\linewidth}
\textbf{Type of isomorphism $\bigM_V \cong \bigM_W$}
\end{minipage}
&
\begin{minipage}[t]{0.21\linewidth}
\textbf{Type of isomorphism $V \cong W$}
\end{minipage}
& \textbf{ $\Rightarrow$ }
& \textbf{ $\Leftarrow$ } &
\begin{minipage}[t]{0.15\linewidth}
\textbf{Reference}
\end{minipage}\\
\hline
& Weak-$*$ continuous  &
\begin{minipage}[t]{0.21\linewidth}
Multiplier biholomorphic
\end{minipage}& $\surd$ &  $\times$  &
\begin{minipage}[t]{0.15\linewidth}
Corollary \ref{Cy:weak_isom}\\
Example \ref{Ex:2BlaschkeSeq}
\end{minipage}\\
\hline
&
\begin{minipage}[t]{0.21\linewidth}
Isometric
\end{minipage}
&
\begin{minipage}[t]{0.21\linewidth}
There is $F \in {\rm Aut}(\ball)$ s.t. $F(W)=V$
\end{minipage}
& $\surd$  & $\surd$ &
\begin{minipage}[t]{0.15\linewidth}
Proposition \ref{Tm:Isom=>CIS}\\
Theorem \ref{Tm:CIS}
\end{minipage}\\
\hline
&
\begin{minipage}[t]{0.21\linewidth}
Completely isometric
\end{minipage}
&
\begin{minipage}[t]{0.21\linewidth}
There is $F \in {\rm Aut}(\ball)$ s.t. $F(W)=V$
\end{minipage}
& $\surd$  & $\surd$   &
\begin{minipage}[t]{0.15\linewidth}
Theorem \ref{Tm:CIS}
\end{minipage}\\
\hline
&
\begin{minipage}[t]{0.21\linewidth}
Unitary equivalence
\end{minipage}
&
\begin{minipage}[t]{0.21\linewidth}
There is $F \in {\rm Aut}(\ball)$ s.t. $F(W)=V$
\end{minipage}
& $\surd$  & $\surd$ &
\begin{minipage}[t]{0.15\linewidth}
Theorem \ref{Tm:CIS}
\end{minipage}\\

\hline
\begin{minipage}[t]{0.21\linewidth}
Finite union of irreducible varieties and a discrete variety
\end{minipage} &
Algebraic &
\begin{minipage}[t]{0.21\linewidth}
Multiplier biholomorphic
\end{minipage}& $\surd$ &  $\times$  &
\begin{minipage}[t]{0.15\linewidth}
Theorem \ref{Tm:irr+dis=>bihol}\\
Example \ref{Ex:2BlaschkeSeq}
\end{minipage}\\

\hline
\begin{minipage}[t]{0.21\linewidth}
Irreducible
\end{minipage} &
Algebraic &
\begin{minipage}[t]{0.21\linewidth}
Multiplier biholomorphic
\end{minipage}&
$\surd$ &  ?  &
\begin{minipage}[t]{0.15\linewidth}
Theorem \ref{Tm:irr+dis=>bihol}\\
Subsection \ref{subsec:problem_finite_union}
\end{minipage}\\

\hline
Homogeneous &
Algebraic &
\begin{minipage}[t]{0.21\linewidth}
There is $A \in {\rm GL}_d(\mathbb C)$ s.t. $A(W)=V$
\end{minipage}& $\surd$ &  $\surd$  &
\begin{minipage}[t]{0.15\linewidth}
 Theorem \ref{Tm:HomogeneousSummary}
\end{minipage}\\

\hline
Homogeneous &
Algebraic &
\begin{minipage}[t]{0.21\linewidth}
Biholomorphic
\end{minipage}& $\surd$ &  $\surd$  &
\begin{minipage}[t]{0.15\linewidth}
 Theorem \ref{Tm:HomogeneousSummary}
\end{minipage}\\

\hline
\begin{minipage}[t]{0.21\linewidth}
Images of finite Riemann surfaces under a holomap that extends to be a 1-to-1 $C^2$-map on the boundary
\end{minipage} &
Algebraic &
\begin{minipage}[t]{0.21\linewidth}
Biholomorphism that exteneds to be a 1-to-1 $C^2$-map on the boundary
\end{minipage}& ? &  $\surd$  &
\begin{minipage}[t]{0.15\linewidth}
Corollary \ref{Cy:RiemannIsom}
\end{minipage}\\

\hline
\begin{minipage}[t]{0.21\linewidth}
Embedded discs
\end{minipage} &
Algebraic &
\begin{minipage}[t]{0.21\linewidth}
Biholomorphic
\end{minipage}& $\surd$ &  $\times$  &
\begin{minipage}[t]{0.15\linewidth}
Example \ref{Ex:TwistedDisc}
\end{minipage}\\
\hline\\
\end{tabular}
\caption{Isomorphisms of varieties in $\ball$ for $d<\infty$ corresponding to isomorphisms of the associated multiplier algebras. The first four lines also hold for $d = \infty$ with minor adjustments.}
\end{sidewaystable}
\newpage

\section{Notation and preliminaries}\label{sec:prelim}

\subsection{Basic notation}
It this survey, $d$ always stands for a positive integer or $\infty = \aleph_0$.
The $d$-dimpensional Hilbert space over $\mathbb C$ is denoted by $\mathbb C^d$ (when $d=\infty$, $\mathbb C^d$ stands for $\ell^2$), and $\ball$ denotes the open unit ball of $\mathbb C^d$. When $d=1$, we usually write $\mathbb D$ instead of $\ball$.

\subsection{The Drury-Arveson space}
Let $H^2_d$ be the {\em Drury-Arveson space} (see \cite{Shalit_DruryArveson}).
$H^2_d$ is the reproducing Hilbert space on $\ball$, the unit ball of $\mathbb C^d$, with kernel functions
\[
k_\lambda(z)=\frac 1{1-\langle z, \lambda \rangle} \quad \text{for } z,\lambda \in \ball.
\]
We denote by $\bigM_d$ the multiplier algebra ${\rm Mult}(H^2_d)$ of $H^2_d$.

\subsection{Varieties and their reproducing kernel Hilbert spaces}
We will use the term {\em analytic variety} (or just a {\em variety}) to refer to the common zero set of a family of $H^2_d$-functions.
If $\cE$ is a set of functions on $\ball$ which is contained in $H^2_d$, let
\[
V(\cE):=\{\lambda \in \ball~:~f(\lambda)=0 \text{ for all } f \in \cE\}.
\]
On the contrary, if $S$ is a subset of $\ball$ let
\[
H_S:=\{ f \in H^2_d~:~f(\lambda)=0 \text{ for all } \lambda \in S \},
\]
and
\[
J_S:=\{ f \in \bigM_d~:~f(\lambda)=0 \text{ for all } \lambda \in S \}.
\]

\begin{Pn}[\cite{DRS2}, Proposition 2.1]
Let $\cE$ be a subset of $H^2_d$, and let $V=V(\cE)$. Then
\[
V=V(J_V).
\]
\end{Pn}
Given an analytic variety $V$, we also define
\[
\bigF_V:=\overline{\rm span}\{k_\lambda~:~\lambda \in V\}.
\]
This Hilbert space is naturally a reproducing kernel Hilbert space of functions living on the variety $V$.
\begin{Pn}[\cite{DRS2}, Proposition 2.3]
Let $S\subseteq \ball$. Then
\[
\bigF_S:=\overline{\rm span}\{k_\lambda~:~\lambda \in S\}=\bigF_{V(H_S)}=\bigF_{V(J_S)}.
\]
\end{Pn}

\subsection{The multiplier algebra of a variety}
The reproducing kernel Hilbert
space $\bigF_V$ comes with its multiplier algebra $\bigM_V = {\rm Mult}(\bigF_V )$. This is
the algebra of all functions $f$ on $V$ such that $f h \in \bigF_V$ for all $h \in \bigF_V$.
A standard argument shows that each multiplier determines a bounded
linear operator $M_f \in B(\bigF_V)$ given by $M_fh := fh$.
We will usually identify
the function $f$ with its multiplication operator $M_f$. We will also
identify the subalgebra of $B(\bigF_V)$ consisting of the $M_f$'s and the algebra
of functions $\bigM_V$ (endowed with the same norm).
We let $Z_i$ denote both the multiplier corresponding to the $i$th coordinate function $z \mapsto z_i$, as well as the multiplication operator it gives rise to.
In some cases, for emphasis, we write $Z_i\big|_V$ instead of $Z_i$.

Now consider the map from $\bigM_d$ into $B(\bigF_V)$ sending each multiplier $f$ to $P_{\bigF_V}M_f|_{\bigF_V}$.
One verifies that this map coincides with the map $f \mapsto f|_V$ and therefore its kernel is $J_V$.
Thus, the multiplier norm of $f|_V$, for $f\in \bigM_d$, is $\|f+J_V\|=\|P_{\bigF_V}M_f|_{\bigF_V}\|$. The complete Nevanlinna-Pick property then implies that this map is completely isometric onto $\bigM_V$. This gives rise to the following proposition.

\begin{Pn}[\cite{DRS2}, Proposition 2.6]
Let $V$ be an analytic variety in $\ball$. Then
\[
\bigM_V=\{f|_V~:~f\in\bigM_d\}.
\]
Moreover the mapping $\varphi:\bigM_d \to\bigM_V$ given by $\varphi(f) = f|_V$ induces a
completely isometric isomorphism and weak-$*$ continuous homeomorphism of $\bigM_d/J_V$
onto $\bigM_V$ . For any $g \in \bigM_V$ and any $f\in\bigM_d$ such that $f|_V = g$, we
have $M_g = P_{\bigF_V}M_f |_{\bigF_V}$ . Given any $F\in M_k(\bigM_V)$, one can choose
$\widetilde F \in M_k(\bigM_d)$ so that $\widetilde F|_V = F$ and $\|\widetilde F \|=\|F\|$.
\end{Pn}

In the above proposition we referred to the weak-$*$ topology in $\cM_V$; this is the weak-$*$ topology which $\cM_V$ naturally inherits from $B(\cF_V)$ by virtue of being a \textsc{wot}-closed (hence weak-$*$ closed) subspace.
The fact that $\cM_V$ is a dual space has significant consequences for us.
It is also useful to know the following.

\begin{Pn}[\cite{DRS2}, Lemma 3.1]\label{Pn:weak-*=WOT}
Let $V$ be a variety in $\ball$. Then then weak-$*$ and the weak-operator topologies on $\bigM_V$ coincide.
\end{Pn}

\subsection{The character space of $\bigM_V$}\label{subsec:char}
Let $\cA$ be a unital Banach algebra.
A {\em character} on $\cA$ is a nonzero multiplicative linear functional.
The set of all characters on $\cA$, endowed with the weak-$*$ topology, is called the {\em character space} of $\cA$, and will be denoted by $M(\cA)$.
It is easy to check that a character is automatically unital and continuous with norm $1$.
If furthermore $\cA$ is an operator algebra, then its characters are automatically completely contractive \cite[Proposition 3.8]{Paulsen}.

The algebras we consider are semi-simple commutative Banach algebras, thus one might expect that the maximal ideal space will be a central part of the classification.
However, these algebras are not uniform algebras; moreover, the topological space $M(\cM_V)$ can be rather wild.
Thus the classification does not use $M(\cM_V)$ directly, but rather a subset of characters that can be identified with a subset of $\bB_d$ and can be endowed with additional structure.

Let $V$ be a variety in $\ball$. Since $(Z_1,\dots,Z_d)$ is a row contraction, it holds that
\[
\|(\rho(Z_1),\dots,\rho(Z_d))\| \leq 1 \quad\text{for all }\rho \in M(\bigM_V).
\]
The map $\pi:M(\bigM_V)\to\Cball$, given by
\[
\pi(\rho)=(\rho(Z_1),\dots,\rho(Z_d)),
\]
is continuous as a map from $M(\bigM_V)$, with the weak-$*$ topology, into $\Cball$ (endowed with the weak topology, in case $d=\infty$).
Since $\pi$ is continuous, $\pi(M(\bigM_V))$ is a compact subset of the closed unit ball.
For every $\lambda \in \pi(M(\bigM_V))$, the set $\pi^{-1}\{\lambda\}\subseteq M(\bigM_V)$ is called the {\em fiber} over $\lambda$.

For every $\lambda \in V$, the fiber over $\lambda$ contains the {\em evaluation functional} $\rho_\lambda$, which is given by
\[
\rho_\lambda(f) = f(\lambda) \,\, , \,\, f \in \cM_V.
\]
The following two results are crucial for much of the analysis of the algebras $\cM_V$.
\begin{Pn}[\cite{DRS2}, Proposition 3.2]\label{identification}
$V$ can be identified with the {\sc wot}-continuous characters of $\mathcal M_V$ via the correspondence $\lambda \leftrightarrow \rho_\lambda$.
%
\end{Pn}


\begin{Pn}[\cite{DRS2}, Proposition 3.2]\label{Pn:int}
If $d<\infty$, then
\[
\pi(M(\bigM_V)) \cap \ball=V,
\]
and for every $\lambda \in V$ the fiber over $\lambda$, that is $\pi^{-1}\{\lambda\}$, is a singleton.
\end{Pn}


\subsection{Metric structure in $M(\cM_V)$}
\label{Dn:Pseudohyperbolic}
Let $\nu \in \ball$, and let $\Phi_\nu$ be the automorphism of the ball that exchanges $\nu$ and $0$ (see \cite[p. 25]{Rudin}):
\[
\Phi_\nu(z):=\frac{\nu-P_\nu z - s_\nu Q_\nu z}{1-\langle z,\nu \rangle},
\]
where
\[
P_\nu=
\begin{cases}
\frac{\langle z,\nu\rangle}{\langle \nu,\nu \rangle}\nu & \text{if } \nu \neq 0,\\
0 & \text{if } \nu=0
\end{cases},
\quad
Q_\nu=I-P_\nu,
\,\, \text{ and } \,\,
s_\nu=(1-\|\nu\|^2)^{\frac 12}.
\]
If $\mu \in \bB_d$ is another point, the {\em pseudohyperbolic distance} between $\mu$ and $\nu$ is defined to be
\[
d_{\rm ph}(\mu,\nu):=\|\Phi_\nu(\mu)\|=\|\Phi_\mu(\nu)\|.
\]
One can check that the pseudohyperbolic distance defines a metric on the open ball.

The following proposition will be useful in the sequel.
Among other things it will imply that the metric structure induced on $V$ by the pseudohyperbolic metric is an invariant of $\cM_V$.

\begin{Pn}[\cite{DRS2}, Lemma 5.3]\label{Pn:Pseudohyperbolic}
Let $V$ be a variety in $\ball$.
\begin{enumerate}[(a)]
\item Let $\mu\in \partial\ball$ and let $\varphi \in \pi^{-1}(\mu)$. Suppose that $\psi \in M(\bigM_V)$ satisfies
$\|\psi - \varphi\|<2$. Then $\psi \in \pi^{-1}(\mu)$.

\item If $\mu,\nu \in V$, then
\[
d_{\rm ph}(\mu,\nu)=\frac{\|\rho_\mu-\rho_\nu\|}{\sup_{\|f\|\leq 1}\left|1-f(\mu) \overline{f(\nu)}\right|}.
\]
As a result,
\[
d_{\rm ph}(\mu,\nu) \leq \|\rho_\mu-\rho_\nu\| \leq 2d_{\rm ph}(\mu,\nu).
\]
\end{enumerate}
\end{Pn}


\section{Weak-$*$ continuous isomorphisms}\label{sec:weak}
Let $V$ and $W$ be two varieties in $\ball$.
We say that $V$ and $W$ are {\em biholomorphic} if there exist holomorphic maps $F:\ball \to \mathbb C^{d}$ and  $G:\mathbb B_{d}\to \mathbb C^{d}$ such that $G\circ F|_V={\bf id}_V$ and $F\circ G|_W={\bf id}_W$.
If furthermore the coordinate functions of $F$ are multipliers, then we say that $V$ and $W$ are {\em multiplier biholomorphic}.

In this section we will see that in the finite dimensional case, if there is a weak-$*$ continuous isomorphism between two multiplier algebras $\bigM_V$ and $\bigM_W$, then $V$ and $W$ are multiplier biholomorphic.
We start with the following proposition, which is a basic tool in the theory.

\begin{Pn}[\cite{DRS2}, Proposition 3.4]\label{Pn:DRS2_3.4}
Let $V$ and $W$ be two varieties in $\ball$, and let $\varphi:\bigM_V\to\bigM_W$ be a unital homomorphism. Then $\varphi$ gives rise to a function $F_\varphi:W\to\Cball$ by
\[
F_\varphi=\pi\circ \varphi^*|_W.
\]
Moreover, there exist multipliers $F_1,F_2,\ldots, F_d \in \bigM$ such that
\[
F_\varphi=(F_1|_W, F_2|_W, \ldots, F_d|_W).
\]
Furthermore, if $\varphi$ is completely bounded or $d<\infty$, then $F_\varphi$ extends to a holomorphic function defined on $\ball$.
\end{Pn}
Here and below $\varphi^*$ is the map from $M(\cM_W)$ into $M(\cM_V)$ given by $\varphi^*(\rho) = \rho \circ \varphi$ for all $\rho \in \cM_W$.
\begin{proof}
Proposition \ref{identification} gives rise to the following commuting diagram

\begin{center}
\begin{tikzpicture}
  \node (WotMw) {$\left \{ {\substack{\text{{\sc wot}-continuous}\\ \text{characters of }\mathcal M_W}}\right \}$};
  \node (Mw) [node distance=3.4cm, right of=WotMw] {$M(\mathcal M_W)$};
  \node (Mv) [node distance=3.4cm, right of=Mw] {$M(\mathcal M_V)$};
  \node (WotMv) [node distance=3.4cm, right of=Mv]{$\left \{ {\substack{\text{{\sc wot}-continuous}\\ \text{characters of }\mathcal M_V}}\right \}$};
  \node (W) [below of=WotMw] {$W$};
  \node (Wm) [below of=Mw] {$\pi(M(\bigM_W))$};
  \node (Vm) [below of=Mv] {$\pi(M(\bigM_V))$};
  \node (V) [below of=WotMv] {$V$};
  \draw[right hook->,ultra thick] (WotMw) to node [swap]  {} (Mw);
  \draw[->,ultra thick] (Mw) to node   {$\varphi^*$} (Mv);
  \draw[<-left hook] (Mv) to node [swap]  {} (WotMv);
  \draw[<->,ultra thick] (WotMw) to node [swap]  {\begin{sideways}\tiny $\lambda \leftrightarrow \rho_\lambda$ \end{sideways}} (W);
  \draw[->] (Mw) to node [swap]  {$\pi$} (Wm);
  \draw[->,ultra thick] (Mv) to node [swap]  {$\pi$} (Vm);
  \draw[<->] (WotMv) to node [swap]  {\begin{sideways}\tiny $\lambda \leftrightarrow \rho_\lambda$ \end{sideways}} (V);
  \draw[right hook->] (W) to node [swap]  {} (Wm);
  \draw[<-left hook] (Vm) to node [swap]  {} (V);
\end{tikzpicture}
\end{center}
and the composition of the thick arrows from $W$ to $\overline \pi(M(\bigM_V)) \subseteq \Cball$ yields the map $F_\varphi$.
Now since $\varphi(Z_i)\in \bigM_W=\{f|_W~:~f\in \bigM\}$, there is an element $F_i \in \bigM$ such that $\varphi(Z_i)=F_i|_W$ and $\|F_i\|=\|\varphi(Z_i)\|$.
Thus, for every $\lambda \in W$,
\[
\begin{split}
F_\varphi(\lambda)
&=\pi(\varphi^*(\rho_\lambda))\\
&=\left(\varphi^*(\rho_\lambda)(Z_1), \varphi^*(\rho_\lambda)(Z_2), \ldots, \varphi^*(\rho_\lambda)(Z_d) \right)\\
&=\left(\varphi(Z_1)(\lambda), \varphi(Z_2)(\lambda), \ldots, \varphi(Z_d)(\lambda)\right)\\
&=\left(F_1|_W(\lambda), F_2|_W(\lambda), \ldots, F_d|_W(\lambda)\right).
\end{split}
\]
It remains to show that if $\varphi$ is completely bounded or $d<\infty$ then $(F_1,\dots,F_d)$ defines a function $\ball \to \mathbb C^d$. If $d<\infty$ it is of course clear. If $d=\infty$ and $\varphi$ is completely bounded then
the norm of $\left(\varphi(Z_1), \varphi(Z_2), \ldots\right)$ is finite, and the $F_i$'s could have been chosen such that $\|(M_{F_1},M_{F_2},\ldots)\|=\|\left(\varphi(Z_1), \varphi(Z_2), \ldots\right)\|$. Hence, with this choice of the $F_i$'s, $(F_1,F_2,\dots)$ defines a function $\mathbb B_\infty \to \ell^2$.
\end{proof}
\begin{Rk}\label{Rk:holo}
When $d=\infty$ and $\varphi$ is not completely bounded, we cannot even say that the map $F_\varphi:W\to\Cball$, in the above proposition, is a holomorphic map.
The reason is that by definition a holomorphic function on a variety should be extendable to a holomorphic function on an open neighborhood of the variety. However, it is not clear whether there exists a choice of the $F_i$'s and a neighborhood of $W$ such that for any $\lambda$ in this neighborhood $(F_1(\lambda),F_2(\lambda),\dots)$ belongs to $\ell^2$.
\end{Rk}
Chasing the diagram in the proof of Proposition \ref{Pn:DRS2_3.4} shows that whenever $\varphi^*$ takes weak-$*$ continuous characters of $\bigM_W$ to weak-$*$ continuous characters of $\bigM_V$, $F_\varphi$ maps $W$ into $V$.
Therefore, if $\varphi$ is a weak-$*$ continuous unital homomorphism, then $F_\varphi(W) \subseteq V$.
This, together with the observation that the inverse of a weak-$*$ continuous isomorphism is weak-$*$ continuous, gives rise to the following corollary.
\begin{Cy}[\cite{DRS2}, Corollary 3.6]\label{CY:DRS2_3.6}
Let $V$ and $W$ be varieties in $\ball$. If $\varphi : \bigM_V \to \bigM_W$ is a unital homomorphism that preserves weak-$*$ continuous characters, then $F_\varphi(W) \subseteq V$ and $\varphi$ is given by
\be\label{eq:imp}
\varphi (F)=f\circ F_\varphi  ,\quad f\in \bigM_V.
\ee
Moreover, if there exists a weak-$*$ continuous isomorphism $\varphi:\bigM_V\to\bigM_W$, then $F_\varphi(W)=V$, $F_{\varphi^{-1}}(V)=W$, and there are multipliers $F_1,\dots,F_d , G_1, \dots, G_d\in \bigM$ such that
\[
F_\varphi=(F_1|_W,\dots,F_d|_W), \quad \text{ and } \quad F_{\varphi^{-1}}=(G_1|_V,\dots,G_d|_V).
\]
\end{Cy}
\begin{proof}
It remains only to verify \eqref{eq:imp}, the rest follows from the discussion above.
If $f \in \cM_V$ and $\lambda \in W$, we find
\[
\varphi(f)(\lambda) = \varphi^*(\rho_\lambda)(f) = \rho_{F_\varphi(\lambda)}(f) = f \circ F_\varphi (\lambda),
\]
as required.
\end{proof}
When $d<\infty$, we obtain the following result.
\begin{Cy}[\cite{DRS2}, Corollary 3.8]\label{Cy:weak_isom}
Let $V$ and $W$ be varieties in $\ball$ for $d<\infty$. If there exists a weak-$*$ continuous isomorphism $\varphi:\bigM_V\to\bigM_W$, then $V$ and $W$ are multiplier biholomorphic.
\end{Cy}
The converse does not hold; see Example \ref{Ex:2BlaschkeSeq} (see also Corollary \ref{Cy:isoHinfty}).
We conclude this section with the following assertion which is a direct result of Proposition \ref{Pn:Pseudohyperbolic}(b) together with the fact that isomorphisms are automatically bounded.
\begin{Cy}[\cite{DHS_EmbeddedDiscs}, Theorem 6.2]\label{Cy:BL}
Suppose $F:W \to V$ is a biholomorphism which induces (by composition) an isomorphism $\varphi:\bigM_V\to \bigM_W$.
Then $F$ must be bi-Lipschitz with respect to the pseudohyperbolic metric, i.e., there is a constant $c>0$ such that
\[
c^{-1} d_{\rm ph}(\mu,\nu) \leq d_{\rm ph}(F(\mu),F(\nu)) \leq c d_{\rm ph}(\mu,\nu).
\]
\end{Cy}
The converse does not hold; see \cite[Example 6.6]{DHS_EmbeddedDiscs}.


\section{Isometric, completely isometric, and unitarily implemented isomorphisms}\label{sec:iso}
Let $V$ and $W$ be two varieties in $\ball$. We say that $V$ and $W$ are {\em conformally equivalent} if there exists an automorphism of $\ball$ (that is, a biholomorphism from $\ball$ into itself) which maps $V$ onto $W$.
In this section we will see that if $V$ and $W$ are conformally equivalent then $\bigM_V$ and $\bigM_W$ are (completely) isometrically isomorphic (in fact, unitarily equivalent).
When $d<\infty$ the converse also holds, and morally speaking it also holds for $d=\infty$.
In fact, when $d=\infty$ it may happen that $\cM_V$ and $\cM_W$ are unitarily equivalent but $V$ and $W$ are not conformally equivalent. This, however, can only be the result of an unlucky embedding of $V$ and $W$ into $\bB_\infty$, and is easily fixed.

\subsection{Completely isometric and unitarily implemented isomorphisms}
\begin{Pn}[\cite{DRS2}, Proposition 4.1]\label{Pn:Unitary}
Let $V$ and $W$ be varieties in $\ball$. Let $F$ be an automorphism of $\ball$ that maps $W$ onto $V$.
Then $f \mapsto f \circ F$ is a unitarily implemented completely isometric isomorphism of $\bigM_V$ onto $\bigM_W$; i.e. $M_{f\circ F}=UM_fU^*$. The unitary $U^*$ is the linear extension of the map
\[
U^*k_{w}=c_w k_{F(w)} \quad \text{for } w \in W,
\]
where $c_w=(1-\|F^{-1}(0)\|^2)^{\frac 12} \overline{k_{F^{-1}(0)}(w)}$.
\end{Pn}

The proof in \cite{DRS2} relies on Theorem 9.2 of \cite{DRS1}, which uses Voiculescu's construction of automorphisms of the Cuntz algebra.
For the convenience of the reader we give here a slightly different proof.

\begin{proof}
Let $F$ be such an automorphism, and set $\alpha=F^{-1}(0)$.
We first show that the linear transformation defined on reproducing kernels by $k_w \mapsto c_wk_{F(w)}$ extends to be a bounded operator of norm $1$.
First note that $\overline{c_w^{-1}}=(1-\|\alpha\|^2)^{-\frac 12}(1-\langle w,\alpha \rangle)$, so $\overline{c_w^{-1}}$ (as a function of $w$) is a multiplier.
The transformation formula for ball automorphisms \cite[Theorem 2.2.5]{Rudin}, shows that
\[
k_{F(w)}(F(z))=c_w^{-1}\overline{c_z^{-1}}k_w(z) \quad \text{ for } w,z \in \ball.
\]
Now,
\[
\langle c_w k_{F(w)},c_z k_{F(z)} \rangle
=c_w \ol{c_z} k_{F(w)}(F(z))=k_w(z)=\langle k_w, k_z \rangle.
\]
Thus, the linear transformation  $k_w \mapsto c_wk_{F(w)}$ extends to an isometry. We denote by $U$ its adjoint.
A short calculation shows that
\[
Uh=(1-\|\alpha\|^2)^{\frac 12} k_\alpha \cdot (h \circ F) \quad \text{ for } h\in H^2_d.
\]
We have already noted that $U^*$ is an isometry, and since its range is evidently dense we conclude that $U$ is a unitary.

Finally, we show that conjugation by $U$ implements the isomorphism between $\bigM_V$ and $\bigM_W$ given by composition with $F$. For $f \in \bigM_V$ and $w \in W$,
\[
UM_f^*U^*k_w=UM_f^*c_wk_{F(w)}=\overline{f(F(w))}Uc_w k_{F(w)}=\overline{(f \circ F)(w)}k_w.
\]
Therefore, $M_{f \circ F}$ is a multiplier on $\bigF_W$ and $M_{f \circ F}=UM_fU^*$.\\
\end{proof}

Before discussing the converse direction, we recall a few definitions on affine sets.
The {\em affine span} (or {\em affine hull}) of a set $S \subseteq \mathbb C^d$
is the set ${\rm aff}(S):=\lambda + {\rm span}(S-\lambda)$ for $\lambda \in S$. This is independent of the choice of $\lambda$. An {\em affine set} is a is a set $A$ with $A = {\rm aff}(A)$. The dimension $\dim(A)$ of an affine
set $A$ is the dimension of the subspace $\ol{A-\lambda}$ for $\lambda \in A$, and the codimension ${\rm codim}(A)$ is the dimension of the quotient space $\mathbb C^d/\ol{A-\lambda}$ for $\lambda \in A$. Both definitions, again, are independent of the choice of $\lambda$. By the {\em affine dimension} (resp. {\em codimension}) of a subset $S\subseteq \mathbb C^d$ we mean the dimension (resp. codimension) of ${\rm aff}(S)$. Furthermore, we use the term {\em affine subset of $\ball$} for any intersection $A\cap \ball$, where $A$ is affine in $\mathbb C^d$.
By \cite[Proposition 2.4.2]{Rudin}, automorphisms of the ball map affine subsets of the ball to affine subsets of the ball.
Therefore, we obtain the following lemma.
\begin{La}\label{Aut=>codim}
Let $V$ and $W$ be varieties in $\ball$ and let $F$ be an automorphism of $\ball$ that maps $W$ onto $V$.
Then, $F(\ol{\rm aff}(V)\cap\ball)=\ol{\rm aff}(W) \cap \ball$. In particular, $\ol{\rm aff}(V)$ and $\ol{\rm aff}(W)$ have the same dimension and the same codimension.
\end{La}
\begin{proof}
The first argument is clear, so it suffices to show that an automorphism of the ball preserves dimensions and codimensions of affine subsets. Indeed, as $F$ is a diffeomorphism, its differential at any point of the ball is an invertible linear transformation. Let $A$ be an affine subset of $\ball$ and let $\lambda \in A$. Let $T_\lambda \ball \cong \mathbb C^d$ be the tangent space of $\ball$ at $\lambda$, and let $T_\lambda A \cong A-\lambda$ be the tangent space of $A$ at $\lambda$. As $A$ is a submanifold of $\ball$, we may think of $T_\lambda A$ as a subspace of $T_\lambda \ball$.
Hence, the invertible linear transformation $dF_\lambda$ maps the subspace $T_\lambda A$ onto $T_{F(\lambda)}F(A)$.
We conclude that $T_\lambda A$ and $T_{F(\lambda)}F(A)$ must have the same dimension and the same codimension.
\end{proof}

Proposition \ref{Pn:Unitary} and Lemma \ref{Aut=>codim} imply, in particular, that if there is an automorphism of the ball which sends $W$ onto $V$, then $V$ and $W$ must have the same affine codimension, and this automorphism gives rise to a completely isometric isomorphism of $\bigM_V$ onto $\bigM_W$ (by precomposing this automorphism). The converse is also true: any completely isometric isomorphism of $\bigM_V$ onto $\bigM_W$, for $V$ and $W$ varieties in the ball having the same affine codimension, arises in this way.

\begin{Pn}\label{Pn:CIS=>Aut}
Let $V$ and $W$ be varieties in $\ball$, with the same affine codimension or with $d<\infty$. Then every completely isometric isomorphism $\varphi:\bigM_V \to \bigM_W$ arises as composition $\varphi(f)=f\circ F$ where $F$ is an automorphism of $\ball$ mapping $W$ onto $V$.
\end{Pn}
\begin{proof}
Recall that Proposition \ref{Pn:DRS2_3.4} assures the existence of a holomorphic map $F:\ball \to \overline{\ball}$ representing $\varphi^*\big|_W$.
A deep result of Kennedy and Yang \cite[Corollary 6.4]{KennedyYang} asserts that $\bigM_V$ and $\bigM_W$ have strongly unique preduals.
It then follows that every isometric isomorphism between these algebras, is also a weak-$*$ homeomorphism.
Thus, by Corollary \ref{CY:DRS2_3.6}, $F(W) \subseteq V$ and $\varphi(f) = f \circ F$.
(We note that if $d < \infty$, then we may argue differently: first one shows using the injectivity of $\varphi$  that $F(\ball) \subseteq \ball$, and then one uses the argument $V=\pi(M(\cM_V)) \cap \bB_d$ of Proposition \ref{Pn:int} to obtain that $\varphi$ preserves weak-$*$ continuous characters.)
Similarly, $\varphi^{-1}:\bigM_W \to \bigM_V$ gives rise to a holomorphic map  $G:\ball \to \ball$ such that $G(V) \subseteq W$ and $\varphi^{-1}(g)=g\circ G$. It is clear that $F\circ G|_V={\bf id}|_V$ and $G\circ F|_W={\bf id}|_W$, and so $F(W)=V$.

By Proposition \ref{Pn:Unitary} and Lemma \ref{Aut=>codim}, we may assume that $V$ and $W$ both contain $0$, and that $F(0)=0$.
Some technical several-complex-variables arguments, which we will not present here, now show that $F|_{\ol{\rm span}W \cap \ball}$ is an isometric linear transformation that maps ${\ol{\rm span}W \cap \ball}$ onto ${\ol{\rm span}V \cap \ball}$ (see \cite[Lemma 4.4]{DRS2}).
In particular, $\ol{\rm span}W$ and $\ol{\rm span}V$ have the same dimension. Since they also have the same codimension, we may extend the definition of $F|_{\ol{\rm span}W \cap \ball}$ to a unitary map on $\mathbb C^d$.
This yields the desired automorphism.
\end{proof}

\begin{Rk}\label{Rk:AffCodim}
The original statement of Proposition \ref{Pn:CIS=>Aut} (which appears in \cite[Theorem 4.5]{DRS2}) does not include the requirement that $V$ and $W$ have the same affine codimension.
Example \ref{Ex:Shift} below shows that this requirement is indeed necessary (for the case $d=\infty$). Nonetheless, it is clear that up to an isometric embedding of the original infinite ball  into a ``larger" one, the original statement does hold. For example, if we replace $V$ and $W$ with their images under the embedding $U : (z_1,z_2,\dots) \mapsto (z_1,0,z_2,0,\dots)$, then both $V$ and $W$ have an infinite affine codimension, and it is now true that $\cM_V$ and $\cM_W$ are completely isometrically isomorphic if and only if $V$ and $W$ are conformally equivalent.
\end{Rk}

\begin{Ex}\label{Ex:Shift}
Let $V=\mathbb B_\infty$ and $W=\{(z_1,z_2,z_3,\dots ) \in \mathbb B_\infty: z_1=0\}$. Let
$F:W\to V$ be defined by
\[
F(0,z_2,z_3,\dots )=(z_2,z_3,\dots ).
\]
Then $F$ is a biholomorphism which cannot be extended to an automorphism of $\mathbb B_\infty$. Let $\varphi: \bigM_V \to \bigM_W$ be defined by $\varphi(f)=f\circ F$. Then $\varphi$ is a completely isometric isomorphism of $\bigM_V$ onto $\bigM_W$, which does not arise as a precomposition with an automorphism of the ball. The reason is of course that $V$ has an affine codimension $0$ while $W$ has an affine codimension $1$.
\end{Ex}

Combining Propositions \ref{Pn:Unitary} and \ref{Pn:CIS=>Aut} yields the following result.
\begin{Tm}[\cite{DRS2}, Theorem 4.5]\label{Tm:CIS}
Let $V$ and $W$ be varieties in $\ball$, with the same affine codimension or with $d<\infty$. Then $\bigM_V$ is completely isometrically isomorphic to $\bigM_W$ if and only if there exists an automorphism $F$ of $\ball$ such that $F(W)=V$.
In fact, under these assumptions, every completely isometric isomorphism $\varphi:\bigM_V \to \bigM_W$ arises as composition $\varphi(f)=f\circ F$ where $F$ is such an automorphism. In this case, $\varphi$ is unitarily implemented by the unitary sending the kernel function $k_w \in \bigF_W$ to a scalar multiple of the kernel function $k_{F(w)}\in\bigF_V$.

If $V$ and $W$ are not assumed to have the same affine codimension, then every completely isometric isomorphism $\varphi:\bigM_V \to \bigM_W$ arises as composition with $U^* \circ F \circ U$, where $F \in \Aut(\bB_d)$ and $U$ is the isometry from Remark \ref{Rk:AffCodim}, and is unitarily implemented.
\end{Tm}

\subsection{Isometric isomorphisms}
By Theorem \ref{Tm:CIS} the conformal geometry of $V$ is completely encoded by the operator algebraic structure $\cM_V$ (and vice versa).
It is natural to ask whether the Banach algebraic structure $\cM_V$ also encodes some geometrical aspect of $V$.
It turns out that within the family of irreducible complete Pick algebras,
every isometric isomorphism of $\bigM_V$ and $\bigM_W$ is actually a completely isometric isomorphism, and the results of the previous section apply.

\begin{La}\label{La:PreservePH}
Let $V$ and $W$ be varieties in $\ball$, and suppose that $\varphi:\bigM_V \to \bigM_W$ is an isometric isomorphism. Then $\varphi^*$ maps $W$ onto $V$ and preserves the pseudohyperbolic distance.
\end{La}
\begin{proof}
The first assertion was obtained in the proof of Proposition \ref{Pn:CIS=>Aut}.
It then follows that $\varphi$ is implemented by composition with $\varphi^*\big|_W$.
Using this together with Proposition \ref{Pn:Pseudohyperbolic} (b), one obtains the second assertion.
\end{proof}

The following theorem appears in \cite[Proposition 5.9]{DRS2} with the additional assumption that $d< \infty$. Here we remove this restriction.
\begin{Tm}[\cite{DRS2}, Proposition 5.9]\label{Tm:Isom=>CIS}
Let $V$ and $W$ be varieties in $\ball$.  Then every isometric isomorphism of $\bigM_V$ onto $\bigM_W$ is completely isometric, and thus is unitarily implemented.
\end{Tm}
\begin{proof}
Without the loss of generality we may assume that $V$ and $W$ have the same affine codimension by embedding the original ball in a larger one, if needed (see Remark \ref{Rk:AffCodim}).
Let $\varphi$ be an isometric isomorphism of $\bigM_V$ onto $\bigM_W$.
By Lemma \ref{La:PreservePH}, $\varphi^*$ maps $W$ onto $V$ and preserves the pseudohyperbolic distance.
Let $F=F_{\varphi}$.

As above, we may assume that $0$ belongs to both $V$ and $W$, and that $F(0)=0$. Let $w_1,w_2,\ldots \in W$ be a sequnce spanning a dense subset of $\overline{\rm span}W$.
For every $p\geq 1$ let $v_p=F(w_p)=\varphi^*(w_p)$.
Put $r_p:=\|w_p\|=d_\text{ph}(w_p,0)$. Then $\|v_p\|=d_\text{ph}(v_p,0)=r_p$.
For every $p$ let $h_p(z):=\langle z, \frac {v_p}{r_p} \rangle$.
This is a continuous linear functional (restricted to $V$), and thus lies in $\bigM_V$.
Furthermore, since $(Z_1,Z_2, \dots, Z_d)$ is a row contraction it follows that $\|h_p\|_{\bigM_V}\leq 1$, and so $\|\varphi(h_p)\|_{\bigM_W}\leq 1$.

Now, let $w$ be an arbitrary point in $W$, set $v=F(w) \in V$, and fix $p \geq 1$.
Since, $\varphi(h_p)$ is a multiplier of norm at most $1$ which satisfies
$\varphi(h_p)(0)=0$, $\varphi(h_p)(w_p)=h_p(v_p)$ and $\varphi(h_p)(w)=h_p(v)$,
we have by a standard necessary condition for interpolation \cite[Theorem 5.2]{AglerMcCarthy_PickInterpolation} that
\[
\begin{bmatrix}
1 & 1 & 1\\
1 & 1 &  \frac{1-\overline{\langle v,v_p \rangle}}{1-\langle w_p,w \rangle}  \\
1 &\frac{1-\langle v,v_p \rangle}{1-\langle w,w_p \rangle} &\frac{1-|\langle v,v_p/r_p \rangle|^2}{1-\langle w,w \rangle}
\end{bmatrix} \geq 0.
\]
Examining the determinant we find that $\frac{1-\langle v,v_p \rangle}{1-\langle w,w_p \rangle}=1$. Therefore,
\[
\langle v,v_p \rangle=\langle w,w_p \rangle\quad\text{for all }  p.
\]
In particular, we obtain $\langle v_i,v_j \rangle=\langle w_i,w_j \rangle$ for all $i,j$.
Therefore, there is a unitary operator $U:\overline{\rm span}W \to \overline{\rm span}V$ such that $Uw_i=v_i$ for all $1\leq i \leq k$.
Since ${\rm codim}(\overline{\rm span}W)={\rm codim}(\overline{\rm span}V)$, it can be extended to a unitary operator $U$ on $\mathbb C^d$.
From here one shows that $F$ agrees with the unitary $U$, and hence $\varphi$ is implemented by an automorphism of the ball. Thus, by Proposition \ref{Pn:Unitary}, $\varphi$ is completely isometric and is unitarily implemented.
\end{proof}

\section{Algebraic isomorphisms}\label{sec:alg}
We now turn to study the algebraic isomorphism problem.
It is remarkable that, under reasonable assumptions, purely algebraic isomorphism implies multiplier biholomorphism.
Throughout this section we will assume that $d<\infty$.

\subsection{Varieties which are unions of finitely many irreducible varieties and a discrete variety}\label{subsec:irreducible}
Let $V$ be a variety in the ball. We say that $V$ is  {\em irreducible} if for any regular point $\lambda \in V$, the intersection of zero sets of all multipliers vanishing on a small neighborhood $V\cap B_\epsilon(\lambda)$ is exactly $V$.
We say that $V$ is {\em discrete} if it has no accumulation points in $\ball$.
We will see that if $V$ and $W$ are two varieties in $\ball$ ($d<\infty$), which are the union of finitely many irreducible varieties and a discrete variety, then whenever $\bigM_V$ and $\bigM_W$ are algebraically isomorphic, $V$ and $W$ are multiplier biholomorphic.

\begin{Rk}
The definition of irreducibility given in the previous paragraph is not to be confused with the classical notion of irreducibility (that is, that there is not non-trivial decomposition of the variety into subvarieties).
Nonetheless, whenever a variety $V$ is irreducible in the classical sense, it is also irreducible in our sense (see e.g. \cite[Theorem, H1]{Gunning}).
\end{Rk}

We open this section with two observations. The first is that every homomorphism between multiplier algebras is norm continuous. A general result in the theory of commutative Banach algebras, says that every homomorphism from a Banach algebra into a commutative semi-simple Banach algebra is norm continuous \cite[Proposition 4.2]{Dales_AutomaticContinuity}.
As $\bigM_W$ is easily seen to be semi-simple, it holds that every homomorphism from $\bigM_V$ to $\bigM_W$ is norm continuous.

The second observation relates to isolated characters of a multiplier algebra. Suppose that $\rho$ is an isolated point in $M(\bigM_V)$.
By Shilov's idempotent theorem \cite[Theorem 5]{BonsallDuncan_CompleteNormedAlgebras}, there is a function $0\neq f\in \bigM_V$ such that every character except $\rho$ annihilates $f$.
As $f\neq 0$, there is $\lambda \in V$ such that $f(\lambda)\neq 0$.
And so, $\rho \in \pi^{-1}(V)$. Thus, when $d<\infty$ any isolated character of a multiplier algebra is an evaluation. This gives rise to the following proposition.

\begin{Pn}[\cite{DRS2}, Lemma 5.2]\label{Pn:discrete}
Let $V$ and $W$ be varieties in $\ball$, with $d<\infty$. Let $\varphi:\bigM_V \to \bigM_W$ be an algebra isomorphism. Suppose that $\lambda$ is an isolated point in $W$. Then $\varphi^*(\rho_\lambda)$ is an evaluation functional at an isolated point in $V$.
\end{Pn}

From the first observation above, together with Proposition \ref{Pn:Pseudohyperbolic}, we obtain:
\begin{Pn}\label{Pn:prop}
Let $V$ and $W$ be a varieties in $\ball$, with $d<\infty$, and let $\varphi : \cM_V \rightarrow \cM_W$ be a homomorphism. Let $U$ be a connected subset of $W$.
Then $\varphi^*(\pi^{-1}(U))$ is either a connected subset of $\pi^{-1}(V)$ (with respect to the norm topology induced by $\bigM_V^*$) or contained in a single fiber of the corona
$M(\bigM_V)\setminus\pi^{-1}(V)$.
\end{Pn}

\begin{Pn}[\cite{DRS2}, Corollary 5.4]\label{Pn:Irr=>onto}
Let $V$ and $W$ be varieties in $\ball$, $d<\infty$, and assume that each one is the union of a discrete variety and a finite union of irreducible varieties. Suppose that $\varphi$ is an algebra isomorphism of $\bigM_V$ onto $\bigM_W$.
Then $\varphi^*$ must map $W$ onto $V$.
\end{Pn}
\begin{proof}
Let us write $V = D_V \cup V_1 \cdots \cup V_m$ and $W = D_W \cup W_1 \cup \cdots \cup W_n$, where $D_V$ and $D_W$ are the discrete parts of $V$ and $W$, and $V_i,W_j$ are all irreducible varieties of dimension at least $1$.
By Proposition \ref{Pn:discrete} $\varphi^*$ maps $D_W$ onto $D_V$.

First let us show that if $W_1$, say, is not mapped entirely into $V$ then it is mapped into a single fiber of the corona $M(\bigM_V)\setminus\pi^{-1}(V)$.
Suppose that $\lambda$ is some regular point of $W_1$ mapped to a fiber of the corona.
Without loss of generality, we may assume it is the fiber over $(1,0,...,0)$.
Then the connected component of $\lambda$ in $W_1$ is mapped into the same fiber, by the previous proposition.
If there exists another point $\mu \in W_1$ which is mapped into $V$ or into another fiber in the corona, then by the previous proposition, the whole connected component of $\mu$ is mapped into $V$ or into the other fiber.
The function $h=\varphi(Z_1|_V)-1|_W$ vanishes on the component of $\lambda$, but does not vanish on the component containing $\mu$.
This contradicts the fact that $W_1$ is irreducible.

Thus, to show that $W_1$ is mapped into $V$ we must rule out the possibility that it is mapped into a single fiber of the corona.
Fix $\lambda \in W_1 \setminus \bigcup_{i=2}^n W_i$.
For each $2 \leq i\leq n$, there is a multiplier $h_i \in \bigM_d$ vanishing on $W_i$ and satisfying $h_i(\lambda) \neq 0$.
Moreover, since $D_W$ is a variety, there is a multiplier $k$ vanishing on $D_W$ and satisfying $k(\lambda) \neq 0$.
Hence, $h:=k \prod_{i=2}^n h_i $ belongs to $\bigM_W$ and vanishes on $D_W \cup \bigcup_{i=2}^n W_i$ but not on $W_1$.
Therefore $\varphi^{-1}(h)$ is a non-zero element of $\bigM_V$.

Now suppose that $\varphi^*(W_1)$ is contained in a fiber over a point in $\partial \ball$, say $(1,0,\dots,0)$.
Since $(Z_1-1)|_V$ is never zero, we see that $(Z_1-1)|_V\varphi^{-1}(h)$ is not the zero function.
However, $(Z_1-1)|_V\varphi^{-1}(h)$ vanishes on $\varphi^*(W_1)$. Therefore, $\varphi((Z_1-1)|_V\varphi^{-1}(h))$ vanishes on $W_1$ and on $D_W \cup \bigcup_{i=2}^n W_i$, contradicting the injectivity of $\varphi$. We deduced that $W_1$ is mapped into $V$. Replacing the roles of $V$ and $W$ shows that $\varphi^*$ must map $W$ onto $V$.
\end{proof}
From Proposition and \ref{Pn:Irr=>onto} and Corollary \ref{CY:DRS2_3.6} we obtain the following.

\begin{Tm}[\cite{DRS2}, Theorem 5.6]\label{Tm:irr+dis=>bihol}
Let $V$ and $W$ be varieties in $\ball$, with $d<\infty$, which are each
a union of finitely many irreducible varieties and a discrete variety.
Let $\varphi$ be an algebra isomorphism of $\bigM_V$ onto $\bigM_W$.
Then there
exist holomorphic maps $F$ and $G$ from $\ball$ into $\mathbb C^d$ with coefficients in
$\bigM_d$ such that
\begin{enumerate}[(a)]
\item $F|_W = \varphi^*|_W$ and $G|_V = (\varphi^{-1})^*|_V$,
\item $G \circ F|_W = {\bf id}_W$ and $F \circ G|_V = {\bf id}_V$,
\item $\varphi(f) = f \circ F$ for $f\in \bigM_V$, and
\item $\varphi^{-1}(g)= g \circ G$ for $g \in\bigM_W$.
\end{enumerate}
\end{Tm}
Theorem \ref{Tm:irr+dis=>bihol} shows in particular that every automorphism of $\bigM_d=\bigM_{\ball}$ is implemented as composition by a biholomorphic map of $\ball$ onto itself, i.e. a conformal automorphism of $\ball$. Proposition \ref{Pn:Unitary} shows that these automorphisms are unitarily implemented (hence, completely isometric). Thus, we obtain the following corollary.
\begin{Cy}[\cite{DRS2}, Corollary 5.8]\label{Cr:Md}
Every algebraic automorphism of $\bigM_d$ for $d$ finite is
completely isometric, and is unitarily implemented.
\end{Cy}

The converse of Theorem \ref{Tm:irr+dis=>bihol} does not hold.
\begin{Ex}\label{Ex:2BlaschkeSeq}
Let
\[
V=\left\{1-\frac 1{n^2}: n\in\mathbb N\right\} \quad \text{and} \quad W=\left\{1-e^{-n^2}:n\in \mathbb N\right\}.
\]
Since they both satisfy the Blaschke condition, they are analytic varieties in $\mathbb D$ (recall that $\{a_n \in \mathbb C:n\in\mathbb N\}$ satisfies the Blaschke condition if $\sum\left(1-|a_n|\right)<\infty$). Let $B(z)$ be the Blaschke product with simple zeros at points in $W$. Define
\[
h(z)=1-e^{\frac1{z-1}} \quad\text{and} \quad g(z)=\frac{\log(1-z)+1}{\log(1-z)}\left(1-\frac{B(z)}{B(0)}\right).
\]
Then $g,h \in H^\infty=\bigM_{\mathbb D}$ and they satisfy
\[
h\circ g|_W={\bf id}_W \quad\text{and} \quad g\circ h|_V={\bf id}_V.
\]
However, by the corollary in \cite[p. 204]{Hoffman}, $W$ is an interpolating sequence and $V$ is not. This implies that $\bigM_W$ is algebraically isomorphic to $\ell^\infty$ while $V$ is not (see \cite[Theorem 6.3]{DRS2}). Thus, $\bigM_V$ and $\bigM_W$ cannot be isomorphic.

\end{Ex}

\subsection{Homogeneous varieties}
Let $V$ be a variety in the ball. We say that $V$ is a {\em homogeneous variety} if it is the common vanishing locus of homogeneous polynomials.

We wish to apply Theorem \ref{Tm:irr+dis=>bihol} to homogeneous varieties in $\bB_d$, $d<\infty$.
It is well known that every algebraic variety can be decomposed into a finite union of irreducible varieties, but caution is required, since the well known result is concerned with irreducibility in another sense than the one we used in Section \ref{subsec:irreducible}.
However, one may show that a homogeneous algebraic variety which is irreducible (in the sense of algebraic varieties) is also irreducible in our sense.
\begin{Pn}\label{Pn:homo=>union_of_irred}
Every homogeneous variety in the ball is a union of finitely many irreducible varieties.
\end{Pn}
\begin{proof}
Let $V$ be a homogeneous variety and let $V=V_1\cup \dots \cup V_n$ be its decomposition into {\em algebraic} irreducible homogeneous varieties (in the sense of algebraic varieties).
We will show that every $V_i$ is irreducible in our sense.
By \cite[Theorem E19, Corollary E20]{Gunning}, once we remove the set of singular points $S(V_i)$, the connected components of $V_i\setminus S(V_i)$ is such that their closures are varieties. Since $S(V_i)$ is a homogeneous variety, these connected components are invariant under nonzero scalar multiplication so their closures are homogeneous varieties. Thus, if there was more than one connected component we would obtain an algebraic decomposition of the variety $V_i$, so $V_i \setminus S(V_i)$ is connected. By the identity principle \cite[Theorem, H1]{Gunning}, the $V_i$'s are irreducible in our sense.
\end{proof}
Thus we obtain the following theorem (the original proof of this theorem was somewhat different --- see \cite[Section 11]{DRS2}).
\begin{Tm}[\cite{DRS1}, Theorem 11.7(2)] \label{Tm:homogen_isom=>biholo}
Let $V$ and $W$ be homogeneous varieties in $\ball$, $d<\infty$.
If $\bigM_V$ and $\bigM_W$ are algebraically isomorphic, then there is a multiplier biholomorphism mapping $W$ onto $V$.
\end{Tm}
The rest of this subsection is devoted towards the converse direction.
Remarkably, a stronger result than the converse holds: it turns out that the existence of a biholomorphism from $W$ onto $V$ implies that the algebras are isomorphic.

We will start by showing that whenever a homogeneous variety $W\subseteq \ball$ is the image of homogeneous variety $V \subseteq \ball$ under a biholomorphism, then it is also the image of $V$ under an invertible linear transformation. To see this, we first need to present the notion of the {\em singular nucleus} of a homogeneous variety.
Lemma 4.5 of \cite{DRS1} and its proof say that a homogeneous variety $V$ in $\mathbb C^d$ is either a linear subspace, or has singular points, and that whenever it is not a linear subspace, the set of singular points $S(V)$ (also known as the {\em singular locus}) of $V$ is a homogeneous variety. Since the dimension of $S(V)$ must be strictly less than the dimension of $V$, there exists a smallest integer $n$ such that $S(\dots(S(S(V)))\dots)$ ($n$ times) is empty.
The set
\[
N(V):=\underbrace{S(\dots(S(S(V)))\dots)}_{n-1 \text{ times}}
\]
is called the {\em singular nucleus} of $V$. By the above discussion, it is a subspace of $\mathbb C^d$.
By basic complex differential geometry, a biholomorphism of $V$ onto $W$ must map $N(V)$ onto $N(W)$.

The following lemma --- which seems to be of independent interest --- was used implicitly in \cite{DRS1}, but in fact does not appear anywhere in the literature. The proof follows closely the proof of \cite[Proposition 4.7]{DRS1}.
\begin{La}\label{La:biholo=>0-biholo}
Let $V$ and $W$ be two biholomorphically equivalent homogeneous varieties in $\ball$. Then there exists a biholomorphism $F$ of $V$ onto $W$ that maps $0$ to $0$.
\end{La}

\begin{proof}
Let $G$ be a biholomorphism of $V$ onto $W$. If $N(V)=N(W)=\{0\}$, then $G(0)=0$, and we are done.
Otherwise, $N(V) \cap \ball$ and $N(W) \cap \ball$ are both complex balls of the same dimension, say $d' \leq d$.
As $G$ takes $N(V) \cap \ball$ onto $N(W) \cap \ball$, we may think of $G$ as an automorphism of $\mathbb B_{d'}$.
We can find two discs $D_1 \subseteq N(V)$ and $D_2 \subseteq N(W)$ such that
$G(D_1) = D_2$ (see \cite[Lemma 4.6]{DRS1}). Define
\[
\bigO(0;V):=\{z \in D_1 : z=F(0) \text{ for some automorphism $F$ of $V$} \}
\]
and
\[
\bigO(0;V,W):=\left\{z \in D_2~ :~
\begin{minipage}{0.5\linewidth}
\text{$z=F(0)$ for some biholomorphism}\\
\text{$F$ of $V$ onto $W$}
\end{minipage}
\right\}.
\]
Since homogeneous varieties are invariant under multiplication by complex numbers, it is easy to check that these sets are circular, that is, for every $\mu \in \bigO(0;V)$ and $\nu \in \bigO(0;V,W)$, it holds that $C_{\mu,D_1}:=\{z \in D_1: |z|=|\mu|\} \subseteq \bigO(0;V)$ and $C_{\nu,D_2}:=\{z \in D_2: |z|=|\nu|\} \subseteq \bigO(0;V,W)$.

Now, as $G(0)$ belongs to $\bigO(0;V,W)$, we obtain that $C := C_{G(0),D_2}\subseteq \bigO(0;V,W)$.
Therefore, the circle $G^{-1}(C)$ is a subset of $\bigO(0;V)$. As $\bigO(0;V)$ is circular, every point of the interior of the circle $G^{-1}(C)$ is a subset of $\bigO(0;V)$.
Thus, the interior of the circle $C$ must be a subset of $\bigO(0;V,W)$. We conclude that $0 \in \bigO(0;V,W)$.
\end{proof}

\begin{Pn}\label{Pn:biholo=>linear}
Let $V$ and $W$ be two biholomorphically equivalent homogeneous varieties in $\ball$. Then there is a linear map on $\mathbb C^d$ which maps $V$ onto $W$.
\end{Pn}
\begin{proof}[Sketch of proof]
By Lemma \ref{La:biholo=>0-biholo}, $V$ and $W$ are biholomorphically equivalent via a $0$ preserving biholomorphism; i.e. there exist two holomorphic maps $F$ and $G$ from $\ball$ into $\mathbb C^d$ such that $G\circ F|_V={\bf id}_V$ and $F \circ G|_W = {\bf id}_W$. Cartan's uniqueness theorem says that if there exists a $0$ preserving biholomorphism between two bounded circular {\em regions}, then it must be a restriction of a linear transformation; see \cite[Theorem 2.1.3]{Rudin}. Now, $V$ and $W$ are indeed circular (since they are homogeneous varieties) and bounded, but do not have to be ``regions" (their interior might be empty). Nevertheless, it turns out that adapting the proof of Cartan's uniqueness theorem to the setting of varieties, rather than regions, does work  (see \cite[Theorem 7.4]{DRS1}). Thus, there exists a linear map $A:\mathbb C^d \to \mathbb C^d$ which agrees with $F$ on $V$.
\end{proof}

Up to now we have seen that if $\cM_V$ and $\cM_W$ are isomorphic, then $V$ and $W$ are biholomorphically equivalent; and we have seen that if $V$ and $W$ are biholomorphically equivalent, then there is a linear map sending $V$ onto $W$, and it is not hard to see that this map can be taken to be invertible.
To close the circle, one needs to show that whenever there is an invertible linear transformation mapping a homogeneous variety $W\subseteq \ball$ onto a homogeneous variety $V\subseteq \ball$, we have that $\bigM_V$ and $\bigM_W$ are similar.
In \cite[Section 7]{DRS1}, this statement was proved for a class of varieties which satisfy some extra assumptions (e.g., irreducible varieties, union of two irreducible components, hypersurfaces, and for the case $d\leq 3$).
Later on, in \cite{Hartz} it was shown that these extra assumptions are superfluous, and that the statement holds for all homogeneous varieties.
The main difficulty was in proving the following lemma.

\begin{La}[\cite{Hartz}]\label{La:PrecomposeBounded}
Let $V$ and $W$ be homogeneous varieties in $\ball$, $d<\infty$,
If there is a linear transformation $A: \mathbb C^d\to\mathbb C^d$ that maps $W$ bijectively onto $V$, then the map $C_{A^*}:\bigF_W \to \bigF_V$, given by
\[
C_{A^*}k_\lambda=k_{A\lambda} \quad \text{ for }\lambda\in W,
\]
is a bounded linear transformation from $\bigF_W$ into $\bigF_V$.
\end{La}
We omit the proof of Lemma \ref{La:PrecomposeBounded}.
The crucial step in its proof is to show that whenever $V_1,\dots, V_n$ are subspaces of $\mathbb C^d$, the algebraic sum of the associated Fock spaces
\[
\bigF(V_1)+\dots+\bigF(V_n) \subseteq \bigF(\mathbb C^d)
\]
is closed. In fact, most of \cite{Hartz} is devoted for proving this crucial step.

\begin{Tm}\label{Tm:linear-biholo=>similar}
Let $V$ and $W$ be homogeneous varieties in $\ball$, $d<\infty$.
If there is an invertible linear transformation $A \in {\rm GL}_d(\mathbb C)$ that maps $W$ onto $V$,
then the map $\varphi:\bigM_V\to \bigM_W$, given by
\[
\varphi(f)=f\circ A \quad\text{ for  } f \in \bigM_V,
\]
is a completely bounded isomorphism, and when regarding $\bigM_V$ and $\bigM_W$ as operator algebras acting on $\bigF_V$ and $\bigF_W$, respectively, $\varphi$ is given by
\[
\varphi(M_f)=(C_{A^*})^*M_f(C_{A^*}^{-1})^* \quad \text{ for  } f \in \bigM_V.
\]
Thus, $\bigM_V$ and $\bigM_W$ are similar.
\end{Tm}
\begin{proof}
By Lemma \ref{La:PrecomposeBounded}, both $C_{A^*}$ and $C_{(A^{-1})^*}$ are bounded, and it is clear that $C_{(A^{-1})^*}=(C_{A^*})^{-1}$.
A calculation shows that $M_{f\circ A}=(C_{A^*})^*M_f (C_{A^*}^{-1})^*$.
\end{proof}

We sum up the results of Theorems \ref{Pn:biholo=>linear}, \ref{Tm:homogen_isom=>biholo} and \ref{Tm:linear-biholo=>similar} as follows.
\begin{Tm}[\cite{DRS1,Hartz}]\label{Tm:HomogeneousSummary}
Let $V$ and $W$ be homogeneous varieties in $\ball$ with $d<\infty$.
Then the following are equivalent:
\begin{enumerate}[(a)]
\item $\bigM_V$ and $\bigM_W$ are similar.
\item $\bigM_V$ and $\bigM_W$ are algebraically isomorphic.
\item $V$ and $W$ are biholomorphically equivalent.
\item There is an invertible linear map on $\mathbb C^d$ which maps $W$ onto $V$.
\end{enumerate}
\end{Tm}

If a linear map $A$ maps $V$ onto $W$ this means that $A$ is length preserving on the homogeneous varieties $\tilde{V}$ and $\tilde{W}$, where $\tilde{V}$ is the homogeneous variety such that $V = \tilde{V} \cap \ball$, and likewise $\tilde{W}$.
This does not mean that $A$ is isometric (as Example \ref{Ex:lines} shows), but it is true that $A$ is isometric on the span of every irreducible component of $W$ \cite[Proposition 7.6]{DRS1}.
Combining this fact with Proposition \ref{Pn:Unitary} we obtain the following result, which sharpens Corollary \ref{Cr:Md} substantially.

\begin{Tm}[\cite{DRS1}, Theorem 8.7]
Let $V$ and $W$ be homogeneous varieties in $\bB_d$, $d<\infty$, such that $W$ is either irreducible or a non-linear hypersurface.
If $\cM_V$ and $\cM_W$ are isomorphic, then they are unitarily equivalent.
\end{Tm}

\begin{Ex}\label{Ex:lines}
Suppose that $V$ and $W$ are each given as the union of two (complex) lines.
There is always a linear map mapping $W$ onto $V$ that is length preserving on $W$, thus $\cM_V$ and $\cM_W$ are algebraically isomorphic.
On the other hand, these algebras will be isometrically isomorphic if and only if the angle between the two lines is the same in each variety.

The case of three lines is also illuminating: it reveals how the algebra $\alg(1,Z)$ and its \textsc{wot}-closure, the algebra $\cM_V$, each encodes different geometrical information.
Indeed, suppose that $V = \spn \{v_1\} \cup \spn\{v_2\} \cup \spn\{v_3\}$ and $W = \spn\{w_1\} \cup \spn\{w_2\} \cup \spn \{w_3\}$, where $v_i,w_j$ are all unit vectors in $\bC^2$ spanning distinct lines.
There always exists a bijective linear map from $W$ onto $V$: indeed, define
$$A: w_1 \mapsto a_1 v_1 \, , \, w_2 \mapsto a_2 v_2 ,$$
and choose $a_1, a_2$ so that $w_3 = b_1 w_1 + b_2 w_2$ is mapped to $v_3$.
One only has to choose $a_1,a_2$ such that $a_1b_1v_1 + a_2 b_2 v_2 = v_3$.
It follows that the algebras $\alg(1,Z\big|_V)$ and $\alg(1,Z\big|_W)$ are isomorphic (the latter two algebras are easily seen to be isomorphic to the coordinate rings of the varieties).

On the other hand, if we require the linear map $A$ to be length preserving on $W$, then $|a_1|=|a_2|=1$.
If $v_3 = c_1 v_1 + c_2 v_2$, then for such a map to exist we will need $a_1 b_1 = c_1$ and $a_2 b_2 = c_2$.
This is possible if and only if $|b_1| = |c_1|$ and $|b_2| = |c_2|$.
Thus the algebras $\cM_V$ and $\cM_W$ in this setup are rarely isomorphic.
\end{Ex}

\subsection{Finite Riemann surfaces}
In seeking a the converse of Theorem \ref{Tm:irr+dis=>bihol}, it is natural to restrict attention to certain well behaved classes of varieties.
In the previous subsection it was shown that the converse of Theorem \ref{Tm:irr+dis=>bihol} holds within the class of homogeneous varieties.
In this subsection we concentrate on generic one-dimensional subvarieites of $\bB_d$, $d<\infty$.

A connected {\em finite Riemann surface} $\Sigma$ is a connected open proper subset
of some compact Riemann surface such that the boundary $\partial \Sigma$ is also the
boundary of the closure and is the union of finitely many disjoint simple
closed analytic curves. A general finite Riemann surface is a finite disjoint
union of connected ones.

Let $\Sigma$ be a connected finite Riemann surface and let $a\in\Sigma$ be some base-point. Let $\omega$ be the harmonic measure
with respect to $a$, i.e. the measure on $\partial \Sigma$ with the property that
\[
u(a)=\int_{\partial \Sigma} u(\zeta)d\omega(\zeta)
\]
for every function $u$ that is harmonic on $\Sigma$ and continuous on $\ol\Sigma$.
We denote by $H^2(\Sigma)$ the closure in $L^2(\omega)$ of the space $A(\Sigma):={\rm Hol}(\Sigma)\cap C(\ol\Sigma)$.
In case that $\Sigma$ is not connected we let $H^2(\Sigma)$ be the direct sum of the $H^2$ spaces of the connected components.

The multiplier
algebra of $H^2(\Sigma)$ is $H^\infty(\Sigma)$, the bounded analytic functions on $\Sigma$. Note that the norm
in $H^2(\Sigma)$ depends on the choice of base-point $a$, but the norm in $H^\infty(\Sigma)$
does not, as it is the supremum of the modulus on $\Sigma$; for more details see \cite{AhrenSarason_HpSpaces}.

We say that a proper holomorphic map $G$ from a finite Riemann
surface $\Sigma$ into a bounded open set $U\subseteq \mathbb C^d$ is a {\em holomap} if there is a finite subset
$\Lambda$ of $\Sigma$ with the property that $G$ is non-singular and injective on $\Sigma \setminus \Lambda$.
We say that $G$ is {\em transversal} at the boundary if
\[
\langle DG(\zeta),G(\zeta) \rangle \neq 0 \quad \text{for all } \zeta \in \partial \Sigma.
\]

The first result on this problem \cite{APV03} showed that if $G : \mathbb D \to W$ is a biholomorphic unramified $C^2$-map that is transversal at the boundary, then there is an isomorphism of multiplier
algebras from $\bigM_\mathbb D = H^\infty(\mathbb D)$ to $\bigM_W$ (the assumptions appearing in \cite{APV03} are slightly weaker --- they only required $C^1$ and did not ask for the map to be unramified --- but it seems that one needs a little more; see \cite[p. 1132]{ARS08}).
This was extended to planar domains in \cite[Section 2.3.6]{{ARS08}},
and to finite Riemann surfaces in \cite{KerrMcCarthyShalit}.
Later, it was proved that a holomorphic
$C^1$ embedding of a finite Riemann surface is automatically transversal at the boundary \cite[Theorem 3.3]{DHS_EmbeddedDiscs}.
Combining this automatic transversaility result with \cite[Theorem 4.2]{KerrMcCarthyShalit} we obtain:

\begin{Tm}[\cite{APV03,ARS08,DHS_EmbeddedDiscs,KerrMcCarthyShalit}]\label{Tm:FRS}
Let $\Sigma$ be a finite Riemann surface and $W$ a variety in $\ball$. Let $G:\Sigma \to \ball$ be a holomap that maps $\Sigma$ onto $W$, is $C^2$ up to $\partial \Sigma$, and is one-to-one on $\partial \Sigma$.
Then the map
\[
\alpha:h \mapsto h\circ G \quad\text{ for } h \in \bigF_W
\]
is an isomorphism from $\bigF_W$ onto $H_G^2(\Sigma):=H^2(\Sigma) \cap \{h\circ G : h\in {\rm Hol}(W)\}$.
Consequently, the map $f \mapsto f \circ G$ implements an isomorphism of $\bigM_W$ onto $H_G^\infty(\Sigma):=H^\infty(\Sigma) \cap \{h\circ G : h\in {\rm Hol}(W)\}$.
\end{Tm}

The main idea of the proof goes back to \cite{APV03}. One first shows that $\alpha$, given by
the formula $h\mapsto h\circ G$, is a well defined bounded and invertible map from $\bigF_W$
onto $H^2_G(\Sigma)$, by computing $\alpha^*$ and $\alpha\alpha^*$, and showing that $\alpha\alpha^*$
is an injective Fredholm operator.
The key trick is to break up $\alpha \alpha^*$ as the sum of a Toeplitz operator and a Hilbert-Schmidt operator (see \cite[Theorem 4.2]{KerrMcCarthyShalit} for details).
Being positive and Fredholm, injectivity
implies invertibility, and the first claim in the theorem follows. A straightforward
computation then shows that the asserted isomorphism between $\bigM_W$ and $H_G^\infty(\Sigma)$
is the similarity induced by $\alpha$.

\begin{Cy}\label{Cy:RiemannIsom}
Let $\Sigma$ be a finite Riemann surface, and let $V$ and $W$ be varieties in $\ball$ such that $W=G(\Sigma)$, where $G:\Sigma \to \ball$ is a holomap which is $C^2$ on $\ol\Sigma$ and is one-to-one on $\partial \Sigma$.
Let $F:W\to V$ be a biholomorphism that extends to be $C^2$ and one-to-one on $\overline W$. Then the map $\varphi:\bigM_V\to\bigM_W$, given by
\[
\varphi(f)=f\circ F \quad\text{ for  }f\in \bigM_V,
\]
is an isomorphism.
\end{Cy}

As an application of the above results, we give the following theorem on extension of bounded holomorphic maps from a one dimensional subvariety of the ball to the entire ball (under rather general assumptions).
Such an extension theorem is difficult to prove using complex-analytic techniques, and it is pleasing to obtain it from operator theoretic considerations.

\begin{Cy}[\cite{APV03} and \cite{KerrMcCarthyShalit}, Corollary 4.12]
Let $W$ be as in Theorem \ref{Tm:FRS}.
Then $\cM_W = H^\infty(W)$, and the norms are equivalent.
Consequently, every $h \in H^\infty(W)$ extends to a multiplier in $\cM_d$, and in particular to a bounded holomorphic function on $\bB_d$.
Moreover, there exists a constant $C$ such that for all $h \in H^\infty(W)$, there is an $\tilde{h} \in \cM_d$ such that $\tilde{h}\big|_W = h$ and $\|\tilde{h}\|_\infty \leq \|\tilde{h}\|_{\cM_d} \leq C \|h\|_\infty$.
\end{Cy}

\subsection{A class of counter-examples}
In the last two subsections we saw classes of varieties, for which (well behaved) biholomorphism of the varieties implies isomorphism of the multiplier algebras.
We now turn to exhibiting a class of examples that show that, in general, biholomorphism of the varieties does not imply that the multiplier algebras are isomorphic.
In particular, these  examples show that biholomorphic varieties need not be multiplier biholomorphic.
\begin{Pn}\label{Pn:CrossBdry}
Suppose that $G:\mathbb D \to \ball$ is a proper injective holomorphic map which extends to a differentiable map on $\mathbb D \cup\{-1,1\}$ such that the extension, also denoted by $G$, satisfies $G(1)=G(-1)$. If $V=G(\mathbb D)$ is a variety, then $G^{-1} \not \in \bigM_V$. In particular, the embedding
\[
\bigM_V \to \bigM_\mathbb D=H^\infty, \quad f\mapsto f \circ G
\]
is not surjective.
\end{Pn}
One way to prove this proposition is to observe that such a map $G$ can not be bi-Lipschitz with respect to the pseudohyperbolic metric, and then invoke Corollary \ref{Cy:BL} (see \cite[Remark 6.3]{DHS_EmbeddedDiscs} for details). For an alternative proof, we refer the reader to \cite[Theorem 5.1]{DHS_EmbeddedDiscs}.

\begin{Ex}\label{Ex:TwistedDisc}
Fix $r\in (0,1)$, and let
\[
b(z)=\frac{z-r}{1-rz}.
\]
Note that $b(1)=1$ and $b(-1)=-1$. Define
\[
G(z)=\frac 1{\sqrt{2}}\left(z^2,b(z)^2\right).
\]
%
It is not hard to verify that this map is a biholomorphism satisfying the hypotheses of Proposition \ref{Pn:CrossBdry}.
Therefore, $\bigM_V \subsetneq H^\infty(V)$, and $G^{-1}$ is not a multiplier.
By Corollary \ref{Cy:MobiusImplements} below we obtain that $\bigM_V$ is not isomorphic to $\bigM_\mathbb D=H^\infty$.
\end{Ex}

\section{Embedded discs in $\bB_\infty$}\label{sec:discs}

\subsection{Some general observations}
In this section we will examine multiplier algebras $\cM_V$ where $V = G(\mathbb{D}) \subseteq \ball$ is a biholomorphic image of a disc via a biholomorphism $G : \mathbb{D} \rightarrow \ball$.
The case that interests us most is $d = \infty$.

\begin{Tm}[\cite{DHS_EmbeddedDiscs}, Theorem 2.5]\label{Tm:EmbeddedDiscs}
Let $V$ and $W$ be two varieties in $\ball$,
biholomorphic to a disc via the maps $G_V$ and $G_W$, respectively.
Furthermore, assume that
\begin{enumerate}[(a)]
\item for every $\lambda \in V$, the fiber $\pi^{-1}\{\lambda\}$ is the singleton $\{\rho_\lambda\}$, and
\item $\pi(M(\bigM_V))\cap \ball =V$.
\end{enumerate}
If $\varphi:\bigM_V \to \bigM_W$ is an algebra isomorphism, then $F=F_\varphi|_W$ is a multiplier biholomorphism $F:W \to V$, such that $\varphi(f)=f\circ F$ for all $f\in\bigM_V$.
\end{Tm}
Here $F = F_\varphi$ is the function provided by Proposition \ref{Pn:DRS2_3.4}.
By saying that $F$ is a {\em multiplier biholomorphism} we mean that
(i) $F = (F_1, F_2, \ldots)$ where every $F_i \in \cM_W$, i.e., is a multiplier, and
(ii) $F$ is holomorphic on $W$, in the sense that for every $\lambda \in W$ there is a ball $B_\epsilon(\lambda)$ and a holomorphic function $\tilde{F} : B_\epsilon(\lambda) \rightarrow \bC^d$ such that $F\big|_{B_\epsilon(\lambda) \cap W} = \tilde{F}\big|_{B_\epsilon(\lambda) \cap W}$.
We require slightly different terminology (compared to Section \ref{sec:weak}) because we are dealing with $d = \infty$, and we are not making any complete boundedness assumptions (see Remark \ref{Rk:holo}).
For more details about holomorphic maps in this setting of discs embedded in $\bB_\infty$ see \cite[Section 2]{DHS_EmbeddedDiscs}.
\begin{proof}
We assume that $d = \infty$.
There are two issues here: we need to prove that $F$ is a biholomorphism, and that $F(W) = V$ in the isomorphic case.
For the first issue, let $\alpha=(\alpha_i)_{i=1}^\infty \in \ell^2$. Then
\[
\langle F\circ G_W(z), \alpha \rangle = \sum_{i=1}^\infty \ol{\alpha_i}h_i(z),
\]
where $h_i(z):=F_i\circ G_W(z)$.
As characters are completely contractive, we have
\[
\sum_{i=1}^\infty |h_i(z)|^2 = \|F(G_W(z))\|^2=\|\rho_{G_W(z)}(\varphi (Z|_W))\|^2 \leq \|Z|_W\|^2=1.
\]
Thus, $\sum_{i=1}^\infty \ol{\alpha_i}h_i$ converges uniformly on $W$ since by the Cauchy-Schwartz inequality,
\[
\sum_{n=N}^\infty \left|\ol {\alpha_n} h_n(z)\right| \leq \left(\sum_{n=N}^\infty |\alpha_n|^2 \right)^\frac{1}{2} \xrightarrow{N\to \infty} 0.
\]
Therefore, $\langle F\circ G_W(\cdot), \alpha \rangle$ is holomorphic for all $\alpha$, and it follows that $F$ is holomorphic (see \cite[Section 2]{DHS_EmbeddedDiscs}).

We now show that the injectivity of $\varphi$ implies that $F$ is not constant, and that this implies  $F(W)\subseteq \bB_\infty$.
Suppose that $F$ is the constant function $\lambda$ ($\lambda \in \ol\ball$).
Then for every $i$ we have $\varphi(\lambda_i-Z_i|_V)=\lambda_i-F_i=0$. By the injectivity of $\varphi$,  $Z_i|_V=\lambda_i$, which is impossible as $V$ is not a singleton.
Thus, $F$ is not constant.
If $\mu=F(\lambda)$ lies in $\partial\bB_\infty$ for some $\lambda \in W$, then $\langle F\circ G_W (\cdot), \mu \rangle$ is a holomorphic function into $\ol{\mathbb D}$, which is equal to $1$ at $\lambda$.
The maximum modulus principle would then imply that this function is constant, so this cannot happen.

In view of the previous paragraph, $F(W)\subseteq \bB_\infty$. Since for every $\lambda \in W$, $\varphi^* (\rho_\lambda) \in \pi^{-1}\{F(\lambda)\} \subseteq \pi^{-1}\bB_\infty$, by the assumptions (a) and (b), we conclude that $F$ maps $W$ into $V$, and therefore (by Corollary \ref{CY:DRS2_3.6}) that $\varphi(f)=f\circ F$. In particular, $\varphi$ is weak-$*$ continuous, and so (as $\varphi$ is an isomorphism) $\varphi^{-1}$ is weak-$*$ continuous too.
Thus, both $\varphi^*$ and $(\varphi^{-1})^*$ map point evaluations to point evaluations. We conclude that $F$ is a biholomorphism, mapping $W$ onto $V$.
\end{proof}

\begin{Rk}
We do not know when precisely conditions (a) and (b) in the above theorem hold.
We do not have an example in which they fail.
We do know that if a variety $V$ in $\mathbb B_\infty$ is the intersection of zero sets of a family of polynomials (or more generally, elements in $\cM_\infty$ that are norm limits of polynomials) then (b) holds (see \cite[Proposition 2.8]{DHS_EmbeddedDiscs}).
\end{Rk}

By a familiar result \cite[p. 143]{Hoffman} the automorphisms of $H^\infty$ are the maps $C_\theta(h):=h\circ \theta$ for some M\"obius map $\theta$ (i.e. $\theta(z)=\lambda\left(\frac{z-a}{1-\ol{a}z}\right)$ for $a\in\mathbb D$, and $\lambda\in \partial \mathbb D$).
If $G$ is a biholomorphic map of the disc onto a variety $V$ in $\ball$, then one can transfer the M\"obius maps to conformal automorphisms of $V$ by sending $\theta$ to $G\circ\theta\circ G^{-1}$.
Since this can be reversed, these are precisely the conformal automorphisms of $V$.
We say that $\bigM_V$ is {\em automorphism invariant} if composition with all these conformal maps yields automorphisms of $\bigM_V$.

\begin{Pn}
Let $V$ and $W$ be two varieties in $\ball$,
biholomorphic to a disc via the maps $G_V$ and $G_W$, respectively.
Assume that $V$ satisfies the conditions (a) and (b) of Theorem \ref{Tm:EmbeddedDiscs}. Let $\varphi:\bigM_V \to \bigM_W$ be an algebra isomorphism. Then there is a M\"obius map $\theta$ such that the diagram
\begin{center}
\begin{tikzpicture}
  \node (Mv) {$\bigM_V$};
  \node (Mw) [node distance=3.4cm, right of=Mv] {$\bigM_W$};
  \node (H1) [below of=Mv] {$H^\infty$};
  \node (H2) [below of=Mw] {$H^\infty$};
  \draw[->] (Mv) to node   {$\varphi$} (Mw);
  \draw[->] (Mv) to node [swap]   {$C_{G_V}$} (H1);
  \draw[->] (Mw) to node  {$C_{G_W}$} (H2);
  \draw[->] (H1) to node [swap]   {$C_\theta$} (H2);
\end{tikzpicture}
\end{center}
commutes.
\end{Pn}
The proof follows by Theorem \ref{Tm:EmbeddedDiscs} and the above discussion. We omit the details.

Suppose that the automorphism $\theta$ can be chosen to be the identity, or equivalently, that $C_F$, where $F=G_V\circ G_W^{-1}$, is an isomorphism of $\bigM_V$ onto $\bigM_W$. Then we will say that $\bigM_V$ and $\bigM_W$ are {\em isomorphic via the natural map}.
\begin{Cy}\label{Cy:MobiusImplements}
Let $V$ and $W$ be two varieties in $\ball$,
biholomorphic to a disc via the maps $G_V$ and $G_W$, respectively.
Assume that $V$ satisfies the conditions (a) and (b) of Theorem \ref{Tm:EmbeddedDiscs}.
If $\bigM_V$ or $\bigM_W$ is automorphism invariant, then $\bigM_V$ and $\bigM_W$ are isomorphic if and only if they are isomorphic via the natural map $C_F$, where $F=G_V\circ G_W^{-1}$. In particular, if $\bigM_V$ is isomorphic to $H^\infty$, then $C_{G_V}$ implements the isomorphism.
\end{Cy}

\subsection{A special class of embeddings}\label{subsec:discs}

We now consider a class of embedded discs in $\mathbb B_\infty$.
The principal goal is to exhibit a large class of multiplier biholomorphic discs in $\bB_\infty$ for which we may classify the obtained multiplier algebras.
Though this goal is not obtained fully, we are able to tell when one of these multiplier algebras is isomorphic to $H^\infty := H^\infty(\mathbb D)$.
Moreover, we obtain an uncountable family of embeddings of the disc into $\bB_\infty$ such that all obtained multiplier algebras are mutually non-isomorphic, while the one dimensional varieties associated with them are all multiplier biholomorphic to each other, via a biholomorphism that extends continuously and one-to-one up to the boundary.

Let $(b_n)_{n=1}^\infty$ be an $\ell^2$-sequence of norm $1$ and $b_1\neq 0$. Define $G:\mathbb D \to \mathbb B_\infty$ by
\[
G(z)=(b_1z,b_2z^2,b_3z^3,\dots)\quad \text{ for }z\in\mathbb D.
\]
Then $G:\mathbb D \to G(\mathbb D)\subseteq \mathbb B_\infty$ is a biholomorphism with inverse $b_1^{-1}Z_1|_{G(\mathbb D)}$ and these maps are multipliers.
Moreover, $G(\mathbb D)$ is a variety because the conditions on the sequence $(b_n)$ (namely, that it has norm $1$ and that $b_1\neq 0$) imply that
\[
V:=V(\{b_nz_1^n-b_1^nz_n~:~n\geq 2\})=G(\mathbb D).
\]
It is easy to see that any two varieties arising this way are multiplier biholomorphic.

\begin{Rk}
One may also consider embeddings similar to the above but with the difference that $\sum|b_n|^2 < 1$, and the results obtained are in some sense analogous to what we describe here, but also contain some surprises.
Since the varieties involved are technically different from those on which we concentrate in this survey, we do not elaborate; the reader is referred to \cite[Section 8]{DHS_EmbeddedDiscs}.
\end{Rk}

Define a kernel on $\mathbb D$ by
\[
k_G(z,w)=\frac{1}{1-\langle G(z), G(w) \rangle} \quad \text{ for } z,w\in\mathbb D,
\]
and let $\bigH_G$ be the Hilbert function space on $\mathbb D$ with reproducing kernel $k_G$. Then we can define a linear map $U:\bigF_V \to \bigH_G$ by $Uh=h\circ G$.
Since
\[
\langle k_{G(z)}, k_{G(w)} \rangle = \frac{1}{1-\langle G(z), G(w) \rangle} = \langle (k_G)_z, (k_G)_w \rangle \quad \text{ for all } z,w \in\mathbb D,
\]
it follows that $Uk_{G(z)}=(k_G)_z$ extends to a unitary map of $\bigF_V$ onto $\bigH_G$.
Hence composition with $G$ determines a unitarily implemented completely isometric isomorphism $C_G:\bigM_V \to {\rm Mult}(\bigH_G)$. Therefore, we can work with multiplier algebras of Hilbert function spaces on the disc rather than the algebras $\bigM_V$ itself.

Now write
\[
k_G(z,w)=\frac{1}{1-\sum_{n=1}^\infty |b_n|^2(z\ol{w})^n}=:\sum_{n=0}^\infty a_n(z\ol{w})^n
\]
for a suitable sequence $(a_n)_{n=0}^\infty$. A direct computation shows that the sequence $(a_n)$ satisfies the recursion
\[
a_0=1\quad\text{and}\quad a_n=\sum_{k=1}^n |b_k|^{2}a_{n-k} \quad \text{for }n\geq 1.
\]
Moreover, $0<a_n\leq 1$ for all $n\in\mathbb N$.

Due to the special form of the kernel $k_G$, we may compute the multiplier norm of monomials in $\bigH_G$.
\begin{La}[\cite{DHS_EmbeddedDiscs}, Lemma 7.2]\label{La:MultNorm}
For every $n\in \mathbb N$, it holds that
\[
\|z^n\|_{{\rm Mult}(\bigH_G)}^2=\|z^n\|_{\bigH_G}^2=\frac 1{a_n}.
\]
\end{La}

We now compare between two varieties embedded discs $V$ and $W$ as above.
We let $(b_n^{\text{\footnotesize{V}}})_{n=1}^\infty$ and $(b_n^W)_{n=1}^\infty$ be two $\ell^2$-sequence of norm $1$ and $b_1^V\neq 0\neq b_1^W$, and define $G_V, G_W :\mathbb D \to \mathbb B_\infty$ by
\[
G_V(z)=(b_1^Vz,b_2^Vz^2,b_3^Vz^3,\dots) \quad \text{and}\quad
G_W(z)=(b_1^Wz,b_2^Wz^2,b_3^Wz^3,\dots).
\]
As before, we consider also the sequences $(a_n^V)_{n=0}^\infty$ and $(a_n^W)_{n=0}^\infty$ which satisfy
\[
k_{G_V}(z,w)=\sum_{n=0}^\infty a_n^V(z\ol{w})^n \quad\text{and}\quad
k_{G_W}(z,w)=\sum_{n=0}^\infty a_n^W(z\ol{w})^n.
\]

\begin{Tm}[\cite{DHS_EmbeddedDiscs}, Proposition 7.5]\label{Tm:isom<=>comparable}
The algebras $\bigM_V$ and $\bigM_W$ are isomorphic via the natural map of composition with $G_V\circ G_W^{-1}$ if and only if the sequences $(a_n^V)$ and $(a_n^W)$ are comparable, i.e., if and only if there is some $c>0$ such that $c^{-1}|a_n^V| \leq |a_n^W| \leq c|a_n^V|$ for all $n$.

Furthermore, if $\pi^{-1}\{\lambda\}=\{\rho_\lambda\}$ for every $\lambda \in W$ and
$\bigM_W$ is automorphism invariant, then $\bigM_V$ and $\bigM_W$ are isomorphic if and only if they are isomorphic via the natural map.
\end{Tm}

\begin{proof}
If $(a_n^V)$ and $(a_n^W)$ are comparable, then by Lemma \ref{La:MultNorm} the norms in $\bigH_{G_V}$ and $\bigH_{G_W}$ of the orthogonal base $\{z^n:n\in\mathbb N\}$ are comparable. Thus, the identity map  is an invertible bounded operator between $\bigH_{G_V}$ and $\bigH_{G_W}$. Therefore, ${\rm Mult}(\bigH_{G_V})={\rm Mult}(\bigH_{G_W})$, so that $\bigM_V$ and $\bigM_W$ are isomorphic via the natural map.

Conversely, if $\bigM_V$ and $\bigM_W$ are isomorphic via the natural map then ${\rm Mult}(\bigH_{G_V})={\rm Mult}(\bigH_{G_W})$. Therefore the identity map is an isomorphism between these two semisimple Banach algebras, so the isomorphism is topological.
By Lemma \ref{La:MultNorm}, the sequences $(a_n^V)$ and $(a_n^W)$ are comparable.

If if $\pi^{-1}\{\lambda\}=\{\rho_\lambda\}$ for every $\lambda \in W$ and
$\bigM_W$ is automorphism invariant, then by Corollary \ref{Cy:MobiusImplements}, this is equivalent to $\bigM_V$ being isomorphic to $\bigM_W$ via any isomorphism.
\end{proof}

Corollary 7.4 of \cite{DHS_EmbeddedDiscs} states that if $\bigM_W$ is automorphism invariant and $\sup_{n\geq 1}({a_n^W}/{a_{n-1}^W})<\infty$, then $\pi^{-1}\{\lambda\}=\{\rho_\lambda\}$ for every $\lambda \in W$. This gives rise to examples in which the second part of Theorem \ref{Tm:isom<=>comparable} is meaningful. For example, the following corollary follows by the above by setting $(b_1^W,b_2^W,b_3^W,\dots )=(1,0,0,\dots)$, and noting that $a_n^W=1$ for all $n\in \mathbb N$.

\begin{Cy}
$\bigM_V$ is isomorphic to $H^\infty$ if and only if the sequence $(a_n^V)$ is bounded below.
\end{Cy}

In terms of the sequence $(b_n)$ the result reads as follows.

\begin{Cy}\label{Cy:isoHinfty}
Let $V=G(\mathbb D)$ where $G(z)=(b_1 z, b_2 z^2, b_3 z^3, \dots)$, where $\|(b_n)\|_{\ell^2}=1$ and $b_1\neq 0$. Then $\bigM_V$ is isomorphic to $H^\infty$ if and only if
\[
\sum_{n=1}^\infty n|b_n|^2 < \infty.
\]
\end{Cy}
\begin{proof}
By the Erd\H{o}s-Feller-Pollard theorem (see \cite[Chapter XIII, Section 11]{Feller}) we know that
\[
\lim_{n\to \infty} a_n =\frac{1}{\sum_{n=1}^\infty n|b_n|^2},
\]
where $1/\infty=0$. Hence, $(a_n)$ is bounded below if and only is the series converges.
\end{proof}

\begin{Ex}[\cite{DHS_EmbeddedDiscs}, Example 7.9]\label{Ex:Hs}
For every $s \in [-1,0]$, consider the reproducing kernel Hilbert spaces $\bigH_s$ with kernel
\[
k^s(z,w)=\sum_{n=0}^\infty (n+1)^s(z\ol{w})^n \quad \text{ for } z,w \in \mathbb D.
\]
It is shown in \cite{DHS_EmbeddedDiscs} that these kernels arise from embeddings as above, and also that these embeddings satisfy all the conditions of Theorem \ref{Tm:isom<=>comparable}.
We have that $a^s_n = (n+1)^s$ in this case, and obviously the sequences $\big((n+1)^s\big)_{n=0}^\infty$ and $\big((n+1)^{s'}\big)_{n=0}^\infty$ are not comparable for $s \neq s'$.
Thus the family of algebras $\mlt(\bigH_s)$ is an uncountable family of multiplier algebras of the type we consider which are pairwise non-isomorphic.
Note that all these algebras live on varieties that are multiplier biholomorphic via a biholomorphism that extends continuously to the boundary.
\end{Ex}


\section{Open problems}\label{sec:open}

Though we have accumulated a body of satisfactory results, and although we have a rich array of examples and counter examples, the isomorphism problem for irreducible Pick algebras is far from being solved.
We close this survey by reviewing some open problems.

\subsection{Finite unions of irreducible varieties} \label{subsec:problem_finite_union}
Theorem \ref{Tm:irr+dis=>bihol} implies that in the case where $V$ and $W$ are finite unions of irreducible varieties in $\ball$ (for $d<\infty$), we have that if $\bigM_V$ and $\bigM_W$ are isomorphic then $V$ and $W$ are multiplier biholomorphic.
It is not known whether the converse holds.
We did see an example of multiplier biholomorphic varieties which are {\em infinite} unions of irreducible varieties but with non-isomorphic multiplier algebras; see Example \ref{Ex:2BlaschkeSeq}.
We also saw an example (Example \ref{Ex:TwistedDisc}) of biholomorphic irreducible varieties, with non-isomorphic multiplier algebras; this, however, was not a {\em multiplier} biholomorphism.
And so the question, whether a multiplier biholomorphism of varieties which are a finite union of irreducible ones implies
that the multiplier algebras are isomorphic, remains unsolved for $d<\infty$ (for $d=\infty$ the answer is {\em no}, see Example \ref{Ex:Hs}).

\subsection{Maximal ideal spaces of multiplier algebras} As we remarked in the introduction, in the case $d=\infty$ there are multiplier algebras $\bigM_V$ for which there are points in $\pi^{-1}\mathbb B_\infty \subseteq M(\bigM_V)$ which are not point evaluations;
similarly, there are also multiplier algebras $\bigM_V$ with characters in fibers over points in $\mathbb B_\infty\setminus V$ \cite[Example 2.4]{DHS_EmbeddedDiscs}.
Nevertheless, when we restrict attention to ``sufficiently nice'' varieties, it might be the case that the characters over the varieties do behave appropriately, in the sense that for every $\lambda \in V$ the fiber $\pi^{-1}\{\lambda\}$ is the singleton $\{\rho_\lambda\}$, and $\pi(M(\bigM_V))\cap \mathbb B_\infty = V$.
In particular, it will be interesting to obtain such a result for the family of discs embedded in $\bB_\infty$ by $G(z) = (b_1 z, b_2 z^2, \ldots)$ as in Section \ref{subsec:discs}.
This will amount to obtaining a better understanding of the maximal ideal space of the algebras $\mlt(\cH_G)$.

\subsection{The correct equivalence relation}
Theorem \ref{Tm:irr+dis=>bihol} says (under some assumptions) that if $\cM_V$ and $\cM_W$ are isomorphic then $V$ and $W$ are multiplier biholomorphic.
We have seen a couple of counter examples showing that the converse is not true, but to clarify the nature of the obstruction let us point out the following: {\em multiplier biholomorphism is not an equivalence relation}, while, on the other hand, isomorphism is an equivalence relation;
see \cite[Remark 6.7]{DHS_EmbeddedDiscs}.
This leads to the problem: describe the equivalence relation $\cong$ on varieties given by
\[
V \cong W \,\textrm{ iff } \, \cM_V \, \textrm{ is isomorphic to } \, \cM_W
\]
in complex geometric terms.

\subsection{Structure theory} The central problem dealt with up to now was the isomorphism problem: when are $\cM_V$ and $\cM_W$ isomorphic (or isometrically isomorphic)?
For isometric isomorphisms the problem is completely resolved: the structure of the Banach algebra $\cM_V$ is completely determined by the conformal structure of $V$.
As for algebraic isomorphisms, we know that the biholomorphic structure of $V$ is an invariant of the algebra $\cM_V$.
This opens the door for a profusion of delicate questions on how to read the (operator) algebraic information from the variety, and vice versa.
For example, how is the dimension of $V$ reflected in $\cM_V$?
If $V$ is a finite Riemann surface with $m$ handles and $n$ boundary components, what in the algebraic structure of $\cM_V$ reflects the $m$ handles and the $n$ boundary components?
What about algebraic-geometric invariants, such as number of irreducible components or degree?

\subsection{Embedding dimension}
A particular question in the flavour of the above broad question, is this: given an irreducible complete Pick algebra $\cA$, what is the minimal $d \in \{1,2,\ldots, \infty\}$ such that $\cA$ is isomorphic to $\cM_V$, with $V \subseteq \bB_d$?
This question is interesting --- and the answer is unknown --- even for the case of the multiplier algebra of the well studied {\em Dirichlet space} $\cD$ (see \cite{ARS11}).

\subsection{Other algebras. Norm closed algebras of multipliers}
The isomorphism problem makes sense on many natural algebras, for examples, one may wonder whether, given two varieties $V,W \subseteq \bB_d$, is it true that the algebra $H^\infty(V)$ is (isometrically) isomorphic to $H^\infty(W)$ precisely when $V$ is biholomorphic to $W$?
Answering this question will require an understanding of the maximal ideal spaces of the bounded analytic functions of a variety.

Another natural class of algebras is given by the norm closures of the polynomials in $\cM_V$,
\[
\cA_V = \ol{\bC[z]}^{\| \cdot \|_{\cM_V}}.
\]
(These algebras are sometimes referred to as the {\em continuous multipliers on} $\cF_V$, but this terminology is misleading since in general $\cA_V \subsetneq C(\ol{V}) \cap \cM_V$; see \cite[Section 5.2]{Shalit_DruryArveson}).
In fact, the isomorphism problem was studied in \cite{DRS1} first for the algebras $\cA_V$.
It was later realized that the norm closed algebras present some delicate difficulties; see \cite[Section 7]{DRS2}.
In fact, subtleties arise already in the case $d=1$; see \cite[Section 8]{DRS2}.

\subsection{Approximation and Nullstellensatz}
One of the problems in studying the isomorphism problem for the norm closed algebras $\cA_V$ is the following (see \cite[Section 7]{DRS2} for an explanation of how these issues relate).
Denote by $\cA_d$ the norm closed algebra generated by the polynomials in $\cM_d$.
Let $V\subseteq \bB_d$ be a variety, and assume that $d<\infty$, and that $V$ is determined by polynomials.
Consider the following ideals $K_V = \{p \in \bC[z] : p\big|_V = 0\}$, $I_V = \{f \in \cA_d : f\big|_V = 0\}$, and $J_V = \{f \in \cM_d : f\big|_V = 0\}$.
A natural question is whether $I_V$ is the norm closure of $K_V$, and whether $J_V$ is  the \textsc{wot}-closure of $I_V$.
In other words, we know that every $f \in I_V$ is the norm limit of polynomials, but does the fact that $f$ vanishes on $V$ imply that it can be approximated in norm using only polynomials from $K_V$?
Likewise, is every function in $J_V$ the limit of a bounded and pointwise convergent sequence of polynomials in $K_V$ (or functions in $I_V$)?

It is very natrual to conjecture that the answer is {\em yes}, and this was indeed proved for homogeneous ideals; see \cite[Corollary 6.13]{DRS2} (see also \cite[Corollary 2.1.31]{RamseyThesis} for the \textsc{wot} case).
As may be expected, this approximation result is closely related to an analytic Nullstellensatz: $\sqrt{\cI} = I (V(\cI))$ (here $\cI$ is some norm closed ideal in $\cA_d$, $V(\cI)$ is the the zero locus of the ideal $\cI$, $I(V(\cI))$ is the ideal of all functions in $\cA_d$ vanishing on $V(\cI)$, and $\sqrt{\cI}$ is an appropriately defined radical; see \cite[Theorem 6.12]{DRS2} and \cite[2.1.30]{RamseyThesis}).
However, we understand very little about these issues in the non-homogeneous case.

\bibliographystyle{plain}

\begin{thebibliography}{10}

\bibitem{AM00} J. Agler and J.E. McCarthy.
\newblock {Complete Nevanlinna--Pick Kernels}.
\newblock {\em J. Funct.\ Anal.}, 175:11--124, 2000.


\bibitem{AglerMcCarthy_PickInterpolation}
J. Agler and J.E. McCarthy.
\newblock {\em {Pick interpolation and Hilbert function spaces}},
\newblock Graduate Studies in Mathematics, 44, Providence,
\newblock RI: American Mathematical Society, 2002.

\bibitem{AhrenSarason_HpSpaces}
P.R.~{Ahern} and D.~{Sarason}.
\newblock {The $H^p$ spaces of a class of function algebras}.
\newblock {\em {Acta Math. }}
117, 123--163, 1967.

\bibitem{APV03}
D. Alpay, M. Putinar and V. Vinnikov.
\newblock {A Hilbert space approach to bounded analytic extension in the ball}.
\newblock {\em {Commun. Pure Appl. Anal.}}, 2(2): 139--145, 2003.

\bibitem{ARS08}
N. Arcozzi, R. Rochberg and E.T. Sawyer.
\newblock {Carleson measures for the Drury-Arveson Hardy space and other Besov-Sobolev spaces on complex balls}.
\newblock {\em {Adv. Math.}}, 218(4): 1107--1180, 2008.

\bibitem{ARS11} N. Arcozzi, R. Rochberg, E.T. Sawyer, and B. Wick.
\newblock {The Dirichlet space: a survey}.
\newblock {\em New York J. Math}, 17:45--86, 2011.


\bibitem{Bers} L. Bers.
\newblock {On rings of analytic functions}.
\newblock {\em Bull. Amer. Math. Soc.}, 54:311--315,
1948.

\bibitem{BonsallDuncan_CompleteNormedAlgebras}
F.F.~{Bonsall} and J.~{Duncan}.
\newblock {\em Complete normed algebras}.
\newblock Ergebnisse der Mathematik und ihrer Grenzgebiete, Band 80. Springer-Verlag, New York-Heidelberg, 1973.



\bibitem{Dales_AutomaticContinuity}
H.G.~{Dales}.
\newblock {Automatic continuity: a survey}.
\newblock {\em {Bull. London Math. Soc.}},
10(2):192--183, 1978.

\bibitem{DHS_EmbeddedDiscs}
K.R.~{Davidson}, M.~{Hartz} and O.M.~{Shalit}.
\newblock {Multipliers of embedded discs}.
\newblock {\em arXiv:1307.3204 [math.OA]}, 2014.
Also in
\newblock {\em {Complex Anal. Oper. Theory}},
to appear.

\bibitem{DHSErr}
K.R.~{Davidson}, M.~{Hartz} and O.M.~{Shalit}.
\newblock {Erratum to: Multpliers of embedded discs}.
\newblock {\em {Complex Anal. Oper. Theory}},
to appear.

\bibitem{DavidsonPitts_ncToeplitz}
K.R. Davidson and D.R. Pitts.
\newblock {The algebraic structure of non-commutative analytic Toeplitz algebras}.
\newblock {\em {Math. Ann.}},
311(2):275--303, 1998.

\bibitem{DavPittsPick} K.R. Davidson and D.R. Pitts.
\newblock {Nevanlinna-Pick interpolation for non-commutative analytic Toeplitz algebras}.
\newblock {\em Integral Equations Operator Theory},
31:321--337, 1998.

\bibitem{DavPitts1} K.R. Davidson and D.R. Pitts.
\newblock {Invariant subspaces and hyper-reflexivity for free semigroup algebras}.
\newblock {\em Proc. Lond. Math. Soc.},
78:401--430, 1999.

\bibitem{DRS1}
K.R.~{Davidson}, C.~{Ramsey} and O.M.~{Shalit}.
\newblock {The isomorphism problem for some universal operator algebras}.
\newblock {\em {Adv. Math.}},
228(1):167--218, 2011.

\bibitem{DRS2}
K.R.~{Davidson}, C.~{Ramsey} and O.M.~{Shalit}.
\newblock {Operator algebras for analytic varieties}.
\newblock {\em arXiv:1201.4072 [math.OA]}, 2014.
Also in
\newblock {\em Trans. Amer. Math. Soc.}, 367:1121--1150, 2015.

\bibitem{Feller}
W.~{Feller}.
\newblock {\em {An introduction to probability theory and its applications}},
Vol. I, 3rd ed.,
Jhon Wiley \& sons Inc., 1968.

\bibitem{GRS}
I.M. Gelfand, D.A. Raikov and G.E. Shilov.
\newblock {Commutative normed rings}.
\newblock {\em Uspehi Matem. Nauk (N. S.)}
1(2):48--146, 1946.

\bibitem{Gunning}
R.C.~{Gunning}.
\newblock{\em {Introduction to holomorphic functions of several variables}}
Vol. II, Local theory.
Monterey, CA, 1990.

\bibitem{Hartz}
M.~{Hartz}.
\newblock {Topological isomorphisms for some universal operator algebras}.
\newblock {\em {J. Funct. Anal.}},
263(11):3564--3587, 2012.

\bibitem{Hoffman}
K.~{Hoffman}.
\newblock {\em {Banach spaces of analytic functions}}.
\newblock Prentice-Hall Series in Modern Analysis Prentice-Hall, Inc., Englewood Cliffs, N. J. 1962.

\bibitem{KennedyYang}
M.~{Kennedy} and D.~{Yang}.
\newblock {A non-self-adjoint Lebesgue decomposition}.
\newblock {\em { Anal. PDE}},
7(2):497--512,
2014.

\bibitem{KerrMcCarthyShalit}
M.~{Kerr}, J.E.~{McCarthy} and O.M.~{Shalit}.
\newblock {On the isomorphism question for complete pick multiplier algebras}.
\newblock {\em {Integral Equations Operator Theory}},
76(1):39--53, 2013.

\bibitem{Paulsen} V.I. Paulsen.
\newblock{\em Completely Bounded Maps and Operator Algebras}.
Cambridge University Press, Cambridge, 2002.

\bibitem{Pick}
G. Pick.
\newblock{ Uber die Beschrankungen analytischer Funktionen, welche durch vorgegebene Funktionswerte bewirkt werden}.
\newblock {\em Math Ann.}, 77(1):7--23, 1915.

\bibitem{Pop06} G. Popescu.
\newblock {Operator theory on noncommutative varieties}.
\newblock {\em Indiana Univ.\ Math.\ J.},
5(2):389--442, 2006.

\bibitem{RamseyThesis} C. Ramsey.
\newblock {\em Maximal ideal space techniques in non-selfadjoint operator algebras}.
PhD. thesis, University of Waterloo, 2013. http://hdl.handle.net/10012/7464.

\bibitem{Rudin}
W.~{Rudin}.
\newblock {\em {Function theory in the unit ball of $\mathbb C^n$}}.
\newblock Springer-Verlag, 1980.

\bibitem{Shalit_DruryArveson}
O.M.~{Shalit}.
\newblock {\em {Operator theory and function theory in Drury-Arveson space and its quotients}}.
\newblock  to appear in {\em {Handbook of Operator Theory}},
editor Daniel Alpay, Springer.

\bibitem{ShalitSolel} O.M. Shalit and B. Solel.
\newblock {Subproduct Systems}.
\newblock {\em Doc. Math.},
14:801--868, 2009.

\end{thebibliography}

\def\cprime{$'$} \def\lfhook#1{\setbox0=\hbox{#1}{\ooalign{\hidewidth
  \lower1.5ex\hbox{'}\hidewidth\crcr\unhbox0}}} \def\cprime{$'$}
  \def\cfgrv#1{\ifmmode\setbox7\hbox{$\accent"5E#1$}\else
  \setbox7\hbox{\accent"5E#1}\penalty 10000\relax\fi\raise 1\ht7
  \hbox{\lower1.05ex\hbox to 1\wd7{\hss\accent"12\hss}}\penalty 10000
  \hskip-1\wd7\penalty 10000\box7} \def\cprime{$'$} \def\cprime{$'$}
  \def\cprime{$'$} \def\cprime{$'$}

\end{document}